\documentclass[graybox]{svmult}

\usepackage{mathptmx}       
\usepackage{helvet}         
\usepackage{courier}        
\usepackage{type1cm}        
\usepackage{makeidx}         
\usepackage{graphicx}        
\usepackage{caption}
  \usepackage{amssymb}
  
\usepackage{bm}
\usepackage{booktabs}
\usepackage{amsmath}                             
\usepackage{amssymb}  
\usepackage{multicol}        
\usepackage[bottom]{footmisc}

\makeindex             

\DeclareMathOperator*{\minimize}{minimize}

\begin{document}

\title*{Tensor completion for Color Images and Videos}
\author{Zhen~Long,
	Yipeng~Liu, Longxi~Chen, Ce~Zhu}
\institute{All the authors are with School of Electronic Engineering / Center for Robotics / Center for Information in Medicine, University of Electronic Science and Technology of China (UESTC), Chengdu, 611731, China. (email: yipengliu@uestc.edu.cn).}

%
%
\maketitle

\abstract{Tensor completion recovers missing entries of multiway data. The missing of entries could often be caused during the data acquisition and transformation. In this paper, we provide an overview of recent development in low rank tensor completion for estimating the missing components of visual data, e. g. , color images and videos. First, we categorize these methods into two groups based on the different optimization models. One optimizes factors of tensor decompositions with predefined tensor rank. The other iteratively updates the estimated tensor via minimizing the tensor rank. Besides, we summarize the corresponding algorithms to solve those optimization problems in details. Numerical experiments are given to demonstrate the performance comparison when different methods are applied to color image and video processing.}
\section{Introduction}
\label{sec:1}
Multidimensional data is becoming important for dealing with the science and engineering problems. In these fields, the size of data is increasing. Tensor, which is a generation of matrix, provides a natural way to represent multidimensional data. For example, a color image is a third-order tensor, which has two indexes in the spatial space and one index in the color. Besides, a color video is a forth-order tensor with an additional index in the time variable. Employing tensor analysis to process multidimensional data, such as in signal and image processing~\cite{1,2,3,4,5,6,7,8,9}, Quantum Chemistry~\cite{10,11,12}, computer vision ~\cite{13,14,15} and data mining~\cite{16,17,18}, has aroused more attention in recent years. In most case, some entries of the acquired multidimensional data are missed. Hence, a well-performed recovery technology should be proposed.

Completion tries to recover the missing data from its observed data according to some prior information. With the success of low rank matrix completion~\cite{19,20,21,22,23,24,25,26,27,28}, low rank tensor completion is an extension to process the multi-dimensional data. Owing to the different decomposition formats, a series of tensor completion methods are proposed.
In order to know about the development of the tensor completion in processing color images and videos, we provide a review to cover the recently methods in this field. The related work in~\cite{29} summarizes the tensor completion algorithms in big data analytics. Compared with it, the difference are as follows: (1). Our survey for tensor completion mainly focus on the field in images and videos processing and we do some experiments using the images to test the algorithms; (2). Compared with CP decomposition and Tucker decomposition adopted in previous work, we introduce additional decomposition methods, such as t-SVD, tensor train decomposition and tensor ring decomposition. (3). We mainly divide the tensor completion into two groups. For each group, based on different tensor decomposition methods, we offer several optimization models and algorithms.

The rest of this paper is organized as follows. Section \ref{sec:2} introduces some notations and preliminaries for tensor decomposition. In Section \ref{sec:3}, the matrix completion is reviewed. The detailed overview of tensor completion is presented in Section \ref{sec:4}. Section \ref{sec:5} provides some experiments for tensor completion using different tensor decomposition formats. The conclusion are concluded in Sec. \ref{sec:6}. Finally, the discussion and future work are presented in Sec. \ref{sec:7}.

\section{Notaions and Preliminaries}
\label{sec:2}
 A scalar, a vector, a matrix, and a tensor are written as $ a $, $  \mathbf{a} $,  $  \mathbf{A} $, and $\mathcal{A}$, respectively. $ a^{*} $ denotes the conjugate of $ a $. The matrix nuclear norm of $\mathbf{A}$ is denoted as $\lVert{\bf A}\rVert_{*}=\sum_{n = 1}^N\sigma_{n}(\bf A)$, where $\sigma_{n}(\bf A)$ is the singular value of matrix {\bf A}. For a $ N $-order tensor $ \mathcal{A} \in \mathbb{R}^{ I_{1} \times\cdots\times I_{N}}$, the multi-index is defined as $\overline{i_1i_2\cdots i_N}=i_1+(i_2-1)I_1+(i_3-1)I_1I_2+\cdots+(i_N-1)I_1\cdots I_{N-1}$ for $i_n=1, 2, \cdots, I_n$, $ n=1, 2, \cdots, N $.The vectorization of tensor $\mathcal{A}$ is denoted as $vec({\mathcal{A}})$. The $``\operatorname{reshape}"$ is a reshaping operation from a matrix or a vector to a tensor, e.g., we can obtain a 3rd=order tensor $\mathcal{X}\in \mathbb{R}^{I_{1}\times I_{2}\times I_{3}}$ by the operating $\operatorname{reshape}(\mathbf{X}, I_{1}, I_{2},I_{3})$ from a matrix $\mathbf{X}\in \mathbb{R}^{I_{1}I_{2}\times I_{3}}$ or $\operatorname{reshape}(\mathbf{x}, I_{1}, I_{2},I_{3})$ from a vector $\mathbf{x}\in \mathbb{R}^{I_{1}I_{2}I_{3}}$. The most common typed of tensor operations can be see Table \ref{Tab:1}.
\begin{table}
	\begin{center}
	\caption{Main tensor/matrix products}
	\label{Tab:1}       
	\begin{tabular}{ll}
		\hline
		$\textbf{A}_{<k>} \in \mathbb{R}^{I_1\cdots I_k\times I_{k+1} \cdots I_N}$ & k-unfolding of $\mathcal{A}$   \\
		\hline
		$\textbf{A}_{[k]} \in \mathbb{R}^{I_k\times I_{1} \cdots I_{k-1}I_{k+1}\cdots I_N}$& mode-k unfolding $\mathcal{A}$\\
		\hline
		$c= <\mathcal{A}, \mathcal{B}> $ & inner product of $\mathcal{A}$ and $\mathcal{B}$ with the same size\\
		\hline
		$\mathcal{C}=\mathcal{A}\otimes\mathcal{B}\in \mathbb{R}^{I_1 J_1\times \cdots\times I_N J_N}$& Kronecker product of  $\mathcal{A} \in \mathbb{R}^{I_1\times \cdots\times I_N}$ and $\mathcal{B} \in \mathbb{R}^{J_1\times\cdots\times J_N}$\\
		\hline
		$\mathcal{C}=\mathcal{A}\circ\mathcal{B}\in\mathbb{R}^{I_1\times \cdots\times I_N\times J_1\times\cdots\times J_N}$& tensor or outer product of $\mathcal{A}\in \mathbb{R}^{I_1\times \cdots\times I_N}$ \\&and $\mathcal{B} \in \mathbb{R}^{J_1\times\cdots\times J_M}$\\
		\hline
		$	\mathcal{C}=\mathcal{A}\times_{n}\mathbf{B}\in \mathbb{R}^{I_1\times \cdots\times I_{n-1}\times J\times I_{n+1}\cdots\times I_N}$& mode-n product of  $\mathcal{A}\in \mathbb{R}^{I_1\times \cdots\times I_N}$ and a matrix $\mathbf{B}\in\mathbb{R}^{J\times I_{n}}$ \\
		\hline
		$\mathcal{C}=<\mathcal{A},\mathcal{B}>_L\in \mathbb{R}^{I_1\times \cdots\times I_N\times Q_1\times \cdots\times Q_M}$ & contracted  product of $\mathcal{A} \in \mathbb{R}^{I_1\times \cdots\times I_N\times P_1\times\cdots\times P_L}$ \\
		 &and $\mathcal{B} \in \mathbb{R}^{P_1\times\cdots\times P_L\times Q_1\times \cdots\times Q_M}$\\
		\hline
		$\mathcal{U}=\mathcal{U}_{1}\mathcal{U}_{2}\cdots \mathcal{U}_{N} \in \mathbb{R}^{S_{1}\times (I_1\cdots I_N) \times S_{N+1}}$& connection product of 3-order tensors $\mathcal{U}_{n} \in \mathbb{R}^{S_{n} \times I_{n}\times S_{n+1}}$\\
		\hline
	\end{tabular}
			\end{center}
\end{table}

\begin{definition}
(\textbf{k-unfolding}~\cite{30}) 	For an $N$-order tensor $\mathcal{A} \in \mathbb{R}^{I_{1} \times\cdots\times I_{N}}$ whose k-unfolding is defined as a matrix
\begin{equation}
\textbf{A}_{<k>} \in \mathbb{R}^{I_1\cdots I_k\times I_{k+1} \cdots I_N}\nonumber
\end{equation}
with entries 
\begin{equation}
\textbf{A}_{<k>}(\overline{i_1i_2\cdots i_k},\overline{i_{k+1}\cdots i_N}) = a_{i_1i_2\cdots i_N}\nonumber,
\end{equation}
where $a_{i_1i_2\cdots i_N}$ is the $i_1i_2\cdots i_N$-th element of $\mathcal{A}$.
\end{definition}

\begin{definition}
 (\textbf{Mode-k unfolding}~\cite{30})
 For an  $N$-order tensor $\mathcal{A} \in \mathbb{R}^{I_{1} \times\cdots\times I_{N}}$, its mode-k unfolding is defined as a matrix 
 \begin{equation}
 \textbf{A}_{[k]} \in \mathbb{R}^{I_k\times I_{1} \cdots I_{k-1}I_{k+1}\cdots I_N}\nonumber
 \end{equation}
with entries 
 \begin{equation}
 \textbf{A}_{[k]}( i_k,\overline{i_{1}\cdots i_{k-1}i_{k+1}\cdots i_N}) = a_{i_1i_2\cdots i_N}\nonumber,
 \end{equation}
 	where $a_{i_1i_2\cdots i_N}$ is the $i_1i_2\cdots i_N$-th element of $\mathcal{A}$.	
\end{definition}

\begin{definition}
(\textbf{Tensor inner product})
 For two tensors $\mathcal{A}$ and $\mathcal{B}$ with the same size, their inner product is defined as $<\mathcal{A}, \mathcal{B}> = <\operatorname{vec}(\mathcal{A}), \operatorname{vec}(\mathcal{B})> $,
where $\operatorname{vec}(\mathcal{A})=\operatorname{vec} ({\bf{A}}_{<1>})$.	
\end{definition}

\begin{definition}
(\textbf{Tensor Frobenius norm})
 The Frobenius norm of a tensor is defined as $\Vert\mathcal{A}\Vert_{\operatorname{F}}^2=<\operatorname{vec}(\mathcal{A}),\operatorname{vec}(\mathcal{A})>$.	
\end{definition}

\begin{definition}
 (\textbf{Tensor Kronecker product})
  For multiway arrays $\mathcal{A} \in \mathbb{R}^{I_1\times \cdots\times I_N}$ and $\mathcal{B} \in \mathbb{R}^{J_1\times\cdots\times J_N}$, the Kronecker product can be denoted as $\mathcal{C}=\mathcal{A}\otimes\mathcal{B}\in \mathbb{R}^{I_1 J_1\times \cdots\times I_N J_N}$ with entries $c_{i_{1}j_{1},\dots,i_{N}j_{N}=a_{i_{1},\dots,i_{N}}b_{j_{1},\dots,j_{N}}}$, where $i_{n}j_{n}=j_{n}+(i_{n}-1)J_{n}$.	
\end{definition}

\begin{definition}
 (\textbf{Tensor or outer product})
For two tensors $\mathcal{A}\in \mathbb{R}^{I_1\times \cdots\times I_N}$ and $\mathcal{B} \in \mathbb{R}^{J_1\times\cdots\times J_M}$, the tensor or outer product can be defined as $\mathcal{C}=\mathcal{A}\circ\mathcal{B}\in\mathbb{R}^{I_1\times \cdots\times I_N\times J_1\times\cdots\times J_N}$
 	with entries 
 	$c_{i_{1},\dots,i_{N},j_{1},\dots,j_{M}}=a_{i_{1}\dots,i_{N}}b_{j_{1},\dots,j_{M}}$.	
\end{definition}

\begin{definition}
(\textbf{Tensor mode-n product})
For a tensor $\mathcal{A}\in \mathbb{R}^{I_1\times \cdots\times I_N}$ and a vector $\mathbf{b}\in\mathbb{R}^{I_{n}}$, the tensor mode-n product can be defined as
\begin{equation}
\mathcal{C}=\mathcal{A}\times_{n}\mathbf{b}\in \mathbb{R}^{I_1\times \cdots\times I_{n-1}\times I_{n+1}\cdots\times I_N},\nonumber
\end{equation}
with entries
\begin{equation}
c_{i_{1},\dots,i_{n-1},i_{n+1},\dotsi_{N}}=\sum_{i_{n}=1}^{I_{n}}(a_{i_{1},\dots,a_{N}}b_{i_{n}}).\nonumber	
\end{equation} 
 Meanwhile, for a tensor $\mathcal{A}\in \mathbb{R}^{I_1\times \cdots\times I_N}$ and a matrix $\mathbf{B}\in\mathbb{R}^{J\times I_{n}}$, the tensor mode-n product can be defined as
\begin{equation}
\mathcal{C}=\mathcal{A}\times_{n}\mathbf{B}\in \mathbb{R}^{I_1\times \cdots\times I_{n-1}\times J\times I_{n+1}\cdots\times I_N},\nonumber
\end{equation}
with entries
\begin{equation}
c_{i_{1},\dots,i_{n-1},j,i_{n+1},\dotsi_{N}}=\sum_{i_{n}=1}^{I_{n}}(a_{i_{1},\dots,a_{N}}b_{j,i_{n}}).\nonumber	
\end{equation} 	

This can be expressed in a matrix form $\mathbf{C_{[n]}}=\mathbf{B}\mathbf{A_{[n]}}$ .	
\end{definition}

\begin{definition}
 (\textbf{Tensor contracted  product})
For multiway arrays $\mathcal{A} \in \mathbb{R}^{I_1\times \cdots\times I_N\times P_1\times\cdots\times P_L}$ and $\mathcal{B} \in \mathbb{R}^{P_1\times\cdots\times P_L\times Q_1\times \cdots\times Q_M}$, the contracted product $\mathcal{C} \in \mathbb{R}^{I_1\times \cdots\times I_N\times Q_1\times \cdots\times Q_M}$ is defined as
 \begin{equation}
 \mathcal{C}=<\mathcal{A},\mathcal{B}>_L\nonumber
 \end{equation}
 with entries
 \begin{eqnarray}
 &\quad& c_{i_1,\cdots,i_N,q_1,\cdots,q_M}\nonumber\\
 &=&\sum_{p_1=1,\cdots,p_L=1}^{P_1,\cdots,P_L}a_{i_1,\cdots,i_N,p_1,\cdots,p_L}b_{p_1,\cdots,p_L,q_1,\cdots,q_M}\nonumber.
 \end{eqnarray}	
\end{definition}

\begin{definition}
 (\textbf{CANDECOMP/PARAFAC (CP) decomposition}~\cite{31,32})
 For a multiway array $\mathcal{A} \in \mathbb{R}^{I_{1} \times\cdots\times I_{N}}$, the CP decomposition is defined as
 \begin{eqnarray}
 \mathcal{A}&=\sum_{s=1}^{S} \lambda_{s}\mathbf{b}_{s}^{(1)}\circ\mathbf{b}_{s}^{(2)}\circ\dots\circ\mathbf{b}_{s}^{(N)}\nonumber\\
 &=\mathcal{G}\times_{1}\mathbf{B}^{(1)}\times_{2}\mathbf{B}^{(2)}\dots\times_{N}\mathbf{B}^{(N)},\nonumber
 \end{eqnarray}
where $\mathbf{b}_{s}^{(1)}\circ\mathbf{b}_{s}^{(2)}\circ\dots\circ\mathbf{b}_{s}^{(N)}$ is rank-one tensor, the non-zero entries $\lambda_{s}$ of the diagonal core tensor $\mathcal{G}\in \mathbb{R}^{S\times S\times\dots\times S}$ represent the weight of rank-one tensors and $\mathbf{B}^{(n)}=[\mathbf{b}_{1}^{(n)},\mathbf{b}_{2}^{(n)},\dots, \mathbf{b}_{S}^{(n)}]\in \mathbb{R}^{I_{n}\times S}$. The rank of CP representation is defined as $\operatorname{rank}_{\operatorname{CP}}(\mathcal{A})=S $, which is the sum number of rank-one tensors. For more intuitively understanding CP decomposition, an example for 3-order tensor is shown in Fig. \ref{fig:CP decomposition}.
 \begin{figure}
 	\centering
 	\includegraphics[width=300pt, keepaspectratio]{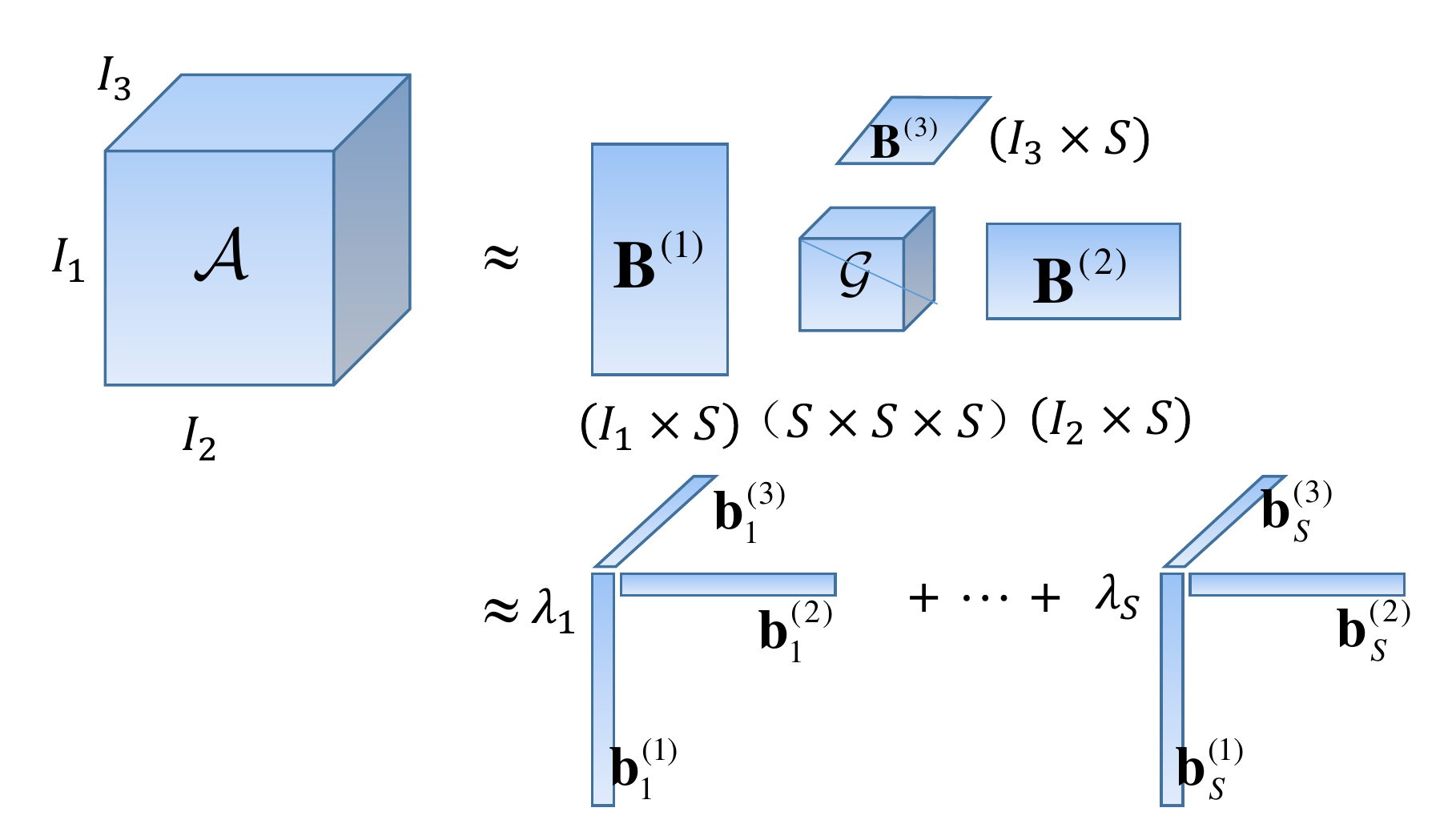}
 	\caption{An example of CP decomposition. }
 	\label{fig:CP decomposition}	
 \end{figure}	
\end{definition}

\begin{definition}
(\textbf{Tucker decomposition}~\cite{33,34,35,32})
For a multiway array $\mathcal{A} \in \mathbb{R}^{I_{1} \times\cdots\times I_{N}}$, the Tucker decomposition is defined as
\begin{eqnarray}
\mathcal{A}&=&\sum_{s_{1}=1}^{S_{1}}\times\sum_{s_{N}=1}^{S_{N}} g_{s_{1}s_{2}\dots s_{N}} (\mathbf{b}_{s_{1}}^{(1)}\circ\mathbf{b}_{s_{2}}^{(1)}\circ\dots\circ\mathbf{b}_{s_{N}}^{(N)})\nonumber\\
&=&\mathcal{G}\times_{1}\mathbf{B}^{(1)}\times_{2}\mathbf{B}^{(2)}\dots\times_{N}\mathbf{B}^{(N)},\nonumber
\end{eqnarray}
where $\mathcal{G}\in \mathbb{R}^{S_{1}\times S_{2}\times\dots\times S_{N}}$, $\mathbf{B}^{(n)}=[\mathbf{b}_{1}^{(n)},\mathbf{b}_{2}^{(n)},\dots, \mathbf{b}_{S}^{(n)}]\in \mathbb{R}^{I_{n}\times S_{n}}$ and the rank of Tucker representation is defined as $\operatorname{rank}_{\operatorname{Tucker}}(\mathcal{A})=(S_{1},S_{2},\cdots, S_{N})$. The detailed illustration of 3-order Tucker decomposition can be seen in Fig. \ref{fig:Tucker decomposition}.
\begin{figure}
	\centering
	\includegraphics[width=250pt, keepaspectratio]{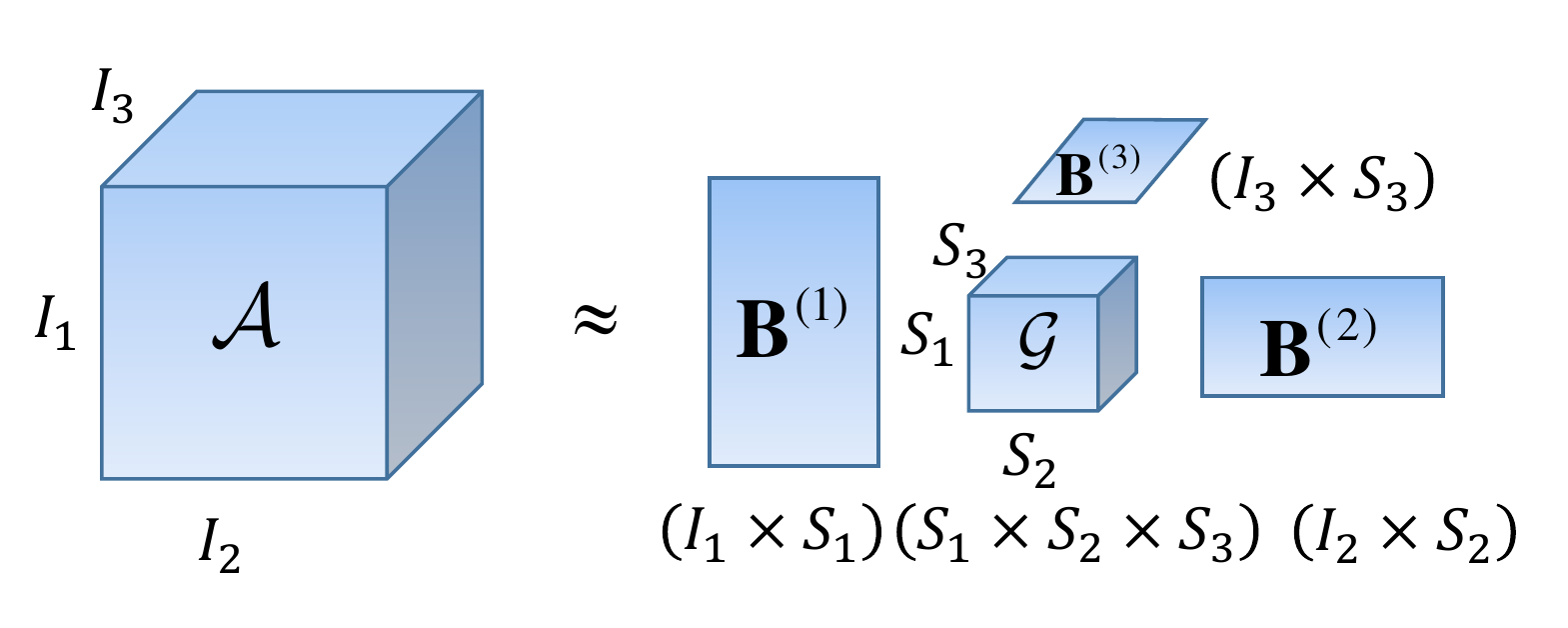}
	\caption{An example of Tucker decomposition. }
	\label{fig:Tucker decomposition}	
\end{figure}		
\end{definition}

\begin{definition}
 (\textbf{Tensor permutation})
For a multiway array $\mathcal{X} \in \mathbb{R}^{I_{1} \times\cdots\times I_{N}}$,  the tensor permutation is defined as $\mathcal{X}^{P_n} \in \mathbb{R}^{I_n\times \cdots\times   I_N\times I_1\cdots\times I_{n-1}}$:
 \begin{equation}
 \mathcal{X}^{P_n}(i_n,\cdots,i_N,i_1,\cdots,i_{n-1})=\mathcal{X}(i_1,\cdots,i_N)\nonumber.
 \end{equation}		
\end{definition}

\begin{definition}
(\textbf{Tensor train decomposition}~\cite{36}) 
For a multiway array $\mathcal{A} \in \mathbb{R}^{I_{1} \times\cdots\times I_{N}}$, the tensor train decomposition is defined as
\begin{equation}
\mathcal{A}(i_1,i_2,\cdots, i_N)=\mathcal{U}_{1}(:,i_1,:)\mathcal{U}_{2}(:,i_2,:)\cdots \mathcal{U}_{N}(:,i_N,:)\nonumber,
\end{equation}	
where the $\mathcal{U}_{n} \in \mathbb{R}^{S_{n}\times I_n\times S_{n+1}} $, $n=1, \cdots, N$ are the cores, and the rank of tensor train representation is defined as $\operatorname{rank}_{\operatorname{TT}}(\mathcal{A})=(S_{1},S_{2},\cdots, S_{N+1})$ with $S_{1}=S_{N+1}=1$. The graphical of tensor train decomposition can be seen in Fig. \ref{fig:TT decomposition}. 
\begin{figure}
	\centering
	\includegraphics[width=300pt, keepaspectratio]{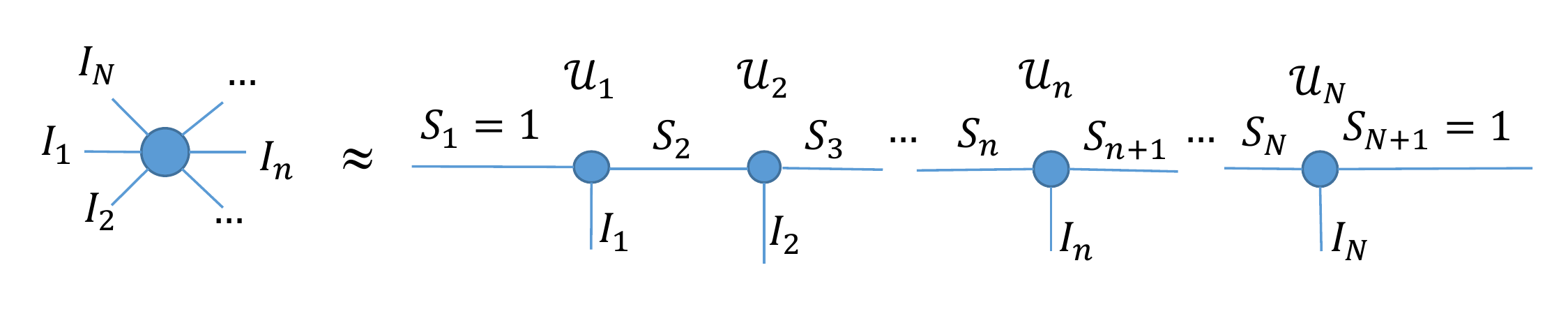}
	\caption{An example of Tensor train decomposition. }
	\label{fig:TT decomposition}	
\end{figure}	
\end{definition}

\begin{definition}
 (\textbf{Tensor ring decomposition}) \cite{37}
 For a multiway array $\mathcal{A} \in \mathbb{R}^{I_{1} \times\cdots\times I_{N}}$, the tensor ring decomposition is defined as
 \begin{equation}
 \mathcal{A}(i_1,i_2,\cdots, i_N)=tr(\mathcal{U}_{1}(:,i_1,:)\mathcal{U}_{2}(:,i_2,:)\cdots \mathcal{U}_{N}(:,i_N,:))\nonumber,
 \end{equation}	
where the $\mathcal{U}_{n} \in \mathbb{R}^{S_{n}\times I_n\times S_{n+1}} $, $n=1, \cdots, N$ are the cores, and the rank of tensor ring representation is defined as $\operatorname{rank}_{\operatorname{TR}}(\mathcal{A})=(S_{1},S_{2},\cdots, S_{N})$ and $S_{N+1}=S_{1}$. The graphical of tensor ring decomposition can be seen in Fig. \ref{fig:TR decomposition}. 
 \begin{figure}
 	\centering
 	\includegraphics[width=250pt, keepaspectratio]{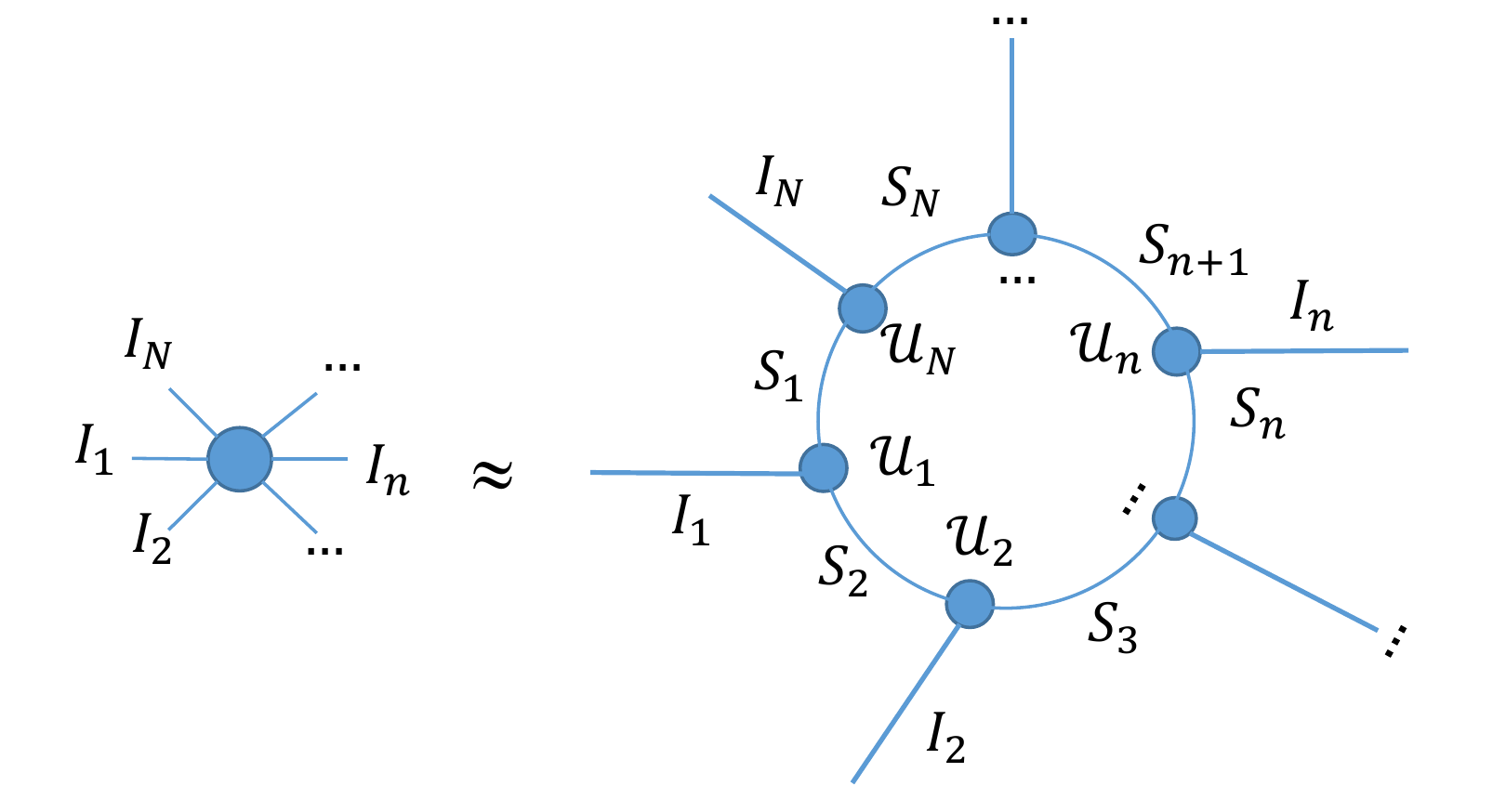}
 	\caption{An example of Tensor ring decomposition. }
 	\label{fig:TR decomposition}	
 \end{figure}
 
 Tensor train is a special case of tensor ring when $S_{1}=1$, we treated these two decomposition as a simple tensor network decomposition format. For an efficient representation of tensor network decomposition, we first introduce a tensor product between 3-order tensors, called tensor connection product.
 
 Before that, the left and right unfoldings of the factors in tensor train decomposition are denoted as:
 \begin{eqnarray}
 {\bf{U}}_{n}^{\operatorname{L}}&=&({\bf{U}}_{n})_{<2>} \in \mathbb{R}^{S_{n}  I_n  \times S_{n+1}}, \nonumber \\
 {\bf{U}}_{n}^{\operatorname{R}}&=&({\bf{U}}_{n})_{<1>} \in \mathbb{R}^{S_{n}\times I_n S_{n+1}}\nonumber.
 \end{eqnarray}		
\end{definition}

\begin{definition}
 (\textbf{Tensor connection product})~\cite{38}
The tensor connection product for $N$ 3-order tensors $\mathcal{U}_{n} \in \mathbb{R}^{S_{n} \times I_{n}\times S_{n+1}}$ is defined as
 \begin{equation}
 \mathcal{U}=\mathcal{U}_{1}\mathcal{U}_{2}\cdots \mathcal{U}_{N} \in \mathbb{R}^{S_{1}\times (I_1\cdots I_N) \times S_{N+1}}\nonumber,
 \end{equation}	
 where $\mathcal{U}_{n}\mathcal{U}_{n+1}  \in \mathbb{R}^{S_{n}\times (I_nI_{n+1}) \times S_{n+2}}$ satisfies
 \begin{equation}
 \mathcal{U}_{n}\mathcal{U}_{n+1}
 = \operatorname{reshape}({\bf{U}}_{n}^{\operatorname{L}}{\bf{U}}_{n+1}^{\operatorname{R}},[S_{n},I_nI_{n+1},S_{n+2}]) \nonumber
 \end{equation}	
 for $n=1,\cdots,N$, where $\operatorname{reshape}$ is a reshaping operation from a matrix of the size $S_{n}I_n\times I_{n+1} S_{n+2}$ to a tensor of the size $S_{n}\times (I_nI_{n+1}) \times S_{n+2}$.
 
 Then simple tensor network decomposition can be represented as $\mathcal{A}=f(\mathcal{U})=f(\mathcal{U}_{1}\mathcal{U}_{2}\cdots \mathcal{U}_{N})$, where function $f$ is a trace operation on $\mathcal{U}(:,i,:)$, $i=1, \cdots, I_1 I_2 \cdots I_N$, and followed by a reshaping operation from vector of the size $1\times(I_1 I_2 \cdots I_N)\times 1$ to tensor of the size $I_1\times I_2\times \cdots\times I_N$.	
\end{definition}

\begin{theorem}
(\textbf{Cyclic permutation property})
Based on the definition of tensor permutation and tensor train decomposition, the tensor permutation of $\mathcal{A}$ is equivalent to the cyclic permutation of its factors in tensor train decomposition form:
\begin{equation*}
\mathcal{A}^{P_n}=f(\mathcal{U}_{n} \cdots\mathcal{U}_{N} \mathcal{U}_{1}\cdots\mathcal{U}_{n-1}),
\end{equation*}
with entries
\begin{eqnarray}
&\mathcal{A}^{P_n}(i_n,\cdots, i_N,i_1,\cdots, i_{n-1})=\operatorname{Trace}(\mathcal{U}_{n}(:,i_n,:) \cdots&\nonumber\\
&\quad\quad \mathcal{U}_{N}(:,i_N,:)\mathcal{U}_{1}(:,i_1,:)\cdots\mathcal{U}_{n-1}(:,i_{n-1},:)).&\nonumber
\end{eqnarray}		
\end{theorem}
\begin{definition}
	(\textbf{t-product}~\cite{39}) For multiway arrays $\mathcal{A}\in \mathbb{R}^{I_{1}\times I_{2}\times I_{3}}$ and $\mathcal{B}\in\mathbb{R}^{I_{2}\times I_{4}\times I_{3}}$, the t-product can be defined as:
	\begin{equation}
	\mathcal{C}=\mathcal{A}\ast\mathcal{B} \in \mathbb{R}^{I_{1}\times I_{4}\times I_{3}} \nonumber
	\end{equation}
	with 
	\begin{equation}
	\mathcal{C}(i_{1}, i_{4}, :)=\sum\limits_{i_{2}=1}^{I_2}\mathcal{A}(i_{1},i_{2},:) \bullet \mathcal{B}(i_{2},i_{4},:),\nonumber
	\end{equation}
	where $\bullet$  denotes the circular convolution between two tubes of same size.
\end{definition}
 \begin{definition}
 	(\textbf{conjugate transpose}) \cite{40} The conjugate transpose of a tensor  $\mathcal{A}$ of size $\emph{I}_1\times\emph{I}_2\times\emph{I}_3$  is the $\emph{I}_2\times\emph{I}_1\times\emph{I}_3$  tensor  $\mathcal{A}^\mathrm{T} $  obtained by conjugate transposing each of the frontal slice and then reversing the order of transposed frontal slices from 2 to $\emph{I}_3$.
 \end{definition}

\begin{definition}
 (\textbf{t-SVD}~\cite{40}) For a tensor $\mathcal{A}\in \mathbb{R}^{I_1\times I_2\times I_3}$ , the t-SVD of  $\mathcal{A}$ is given by
 \begin{equation}
 \mathcal{A}=\mathcal{U}* \mathcal{S}*\mathcal{V^\mathrm{T}}\nonumber
 \end{equation}
 where $\mathcal{U}\in \mathbb{R}^{I_1\times I_1\times I_3}$  and $\mathcal{V}\in\mathbb{R}^{I_2\times I_2\times I_3}$ satisfy $\mathcal{U}\ast\mathcal{U}^{T}=\mathcal{U}^{T}\ast\mathcal{U}=\mathcal{I}\in \mathbb{R}^{I_1\times I_1\times I_3}$and  $\mathcal{V}\ast\mathcal{V}^{T}=\mathcal{V}^{T}\ast\mathcal{V}=\mathcal{I}\in \mathbb{R}^{I_2\times I_2\times I_3}$ respectively, and $\mathcal{S}\in \mathbb{R}^{I_1\times I_2\times I_3}$  is a f-diagonal  tensor whose frontal slices ia a diagonal matrix.
 
 We can obtain this decomposition by computing matrix singular value decomposition (SVD) in the Fourier domain, as it shows in Algorithm 1. Fig. 4 illustrates the decomposition for the three-dimensional case.\smallskip
 
 \begin{table}[htbp]
 	\centering
 	\begin{tabular}{lll}
 		\toprule
 		$\bf Algorithm 1: $ t-SVD for 3-way data \\
 		\midrule
 		\textbf{Input:} $ \mathcal A\in \mathbb{R}^{I_{1}\times I_{2}\times I_{3}} $ \\
 		
 		$\mathcal{D} \gets $fft($\mathcal{A}$,[],$3$),\\
 		
 		$\bf  for$ $ i=1 $ to $ I_{3}$,  \textbf{do}  \\
 		\quad $ [$ \textbf {U , S ,  V }$ ] $ = svd$(\mathcal{D}(:,:,i))$, \\
 		
 		\quad$\hat{\mathcal{U}}(:,:,i)$=$ \bf U $,~ $\hat{\mathcal{S}}(:,:,i)$)=$ \bf S $,~ $\hat{\mathcal{V}}(:,:,i)$ =$\bf V $,\\
 		\bf end for \\
 		
 		$\mathcal{U}\gets$ {ifft($\hat{\mathcal{U}}$,[],3)},~ $\mathcal{S}\gets$ ifft($\hat{\mathcal{S}}$,[],3),~ $\mathcal{V}\gets$ ifft($\hat{\mathcal{V}}$,[],3),\\
 		
 		\textbf{Output:} $\mathcal{U}, \mathcal{S}, \mathcal{V} $\\
 		\bottomrule	
 	\end{tabular}
 \end{table}
 
 \begin{figure}[htbp]
 	\centering
 	\includegraphics[width=260pt, keepaspectratio]{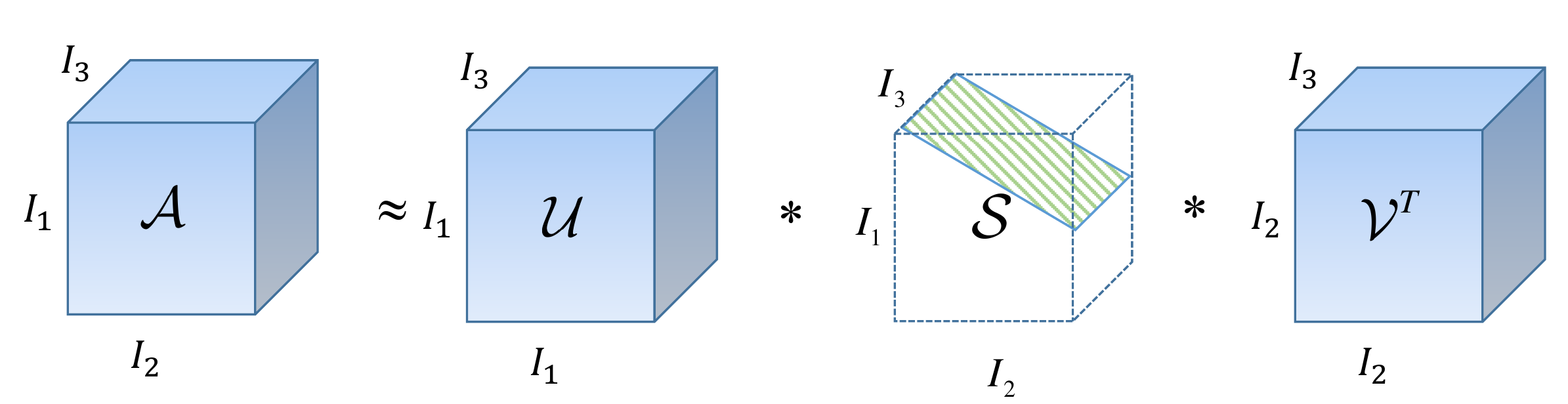}
 	\caption{Illustration of the t-SVD of an $I_1\times I_2\times I_3$ tensor.}
 \end{figure}
	
\end{definition}

\section{Matrix Completion}
\label{sec:3}

Matrix completion is the problem of recovering a data matrix from its partial entries. In many problems, we often assume that the matrix, which we wish to recover, is low rank or approximately low rank. The optimization model for matrix completion was proposed firstly in~\cite{19}, and can be formulated as:
\begin{eqnarray}
&&\minimize_{\mathbf{X}} \quad \operatorname{rank}(\mathbf{X})\nonumber\\
&&\text{subject to}\quad \mathbf{X}_{\mathbb{O}}=\mathbf{T}_{\mathbb{O}} ,
\label{equ:org}
\end{eqnarray}
where $\mathbf{X}$ represents the completed low-rank matrix, the $\operatorname{rank}(\mathbf{X})$ is equal to the rank of the matrix $\mathbf{X}$ and the $\mathbb{O}$ is the entries sets. The equality constraint $\mathbf{X}_{\mathbb{O}}=\mathbf{T}_{\mathbb{O}}$ means that the available entries in observation matrix $\mathbf{T}$ is equal to the available entries in completed matrix $\mathbf{X}$. This optimization problem is NP-hard, and can be solved by some fundamental matrix completion approaches. These methods could be divided into two categories : minimize matrix nuclear norm and low rank matrix decomposition.

\subsection{nuclear norm based matrix completion}
In solving the matrix completion problems, we can use matrix nuclear norm replace non-convex rank function as a convex surrogate. It has been proven that matrix nuclear norm is the tightest lower bound of matrix rank function among all possible convex~\cite{26}. The optimization problem (\ref{equ:org}) can be written as:
\begin{eqnarray}
&& \minimize_{\mathbf{X}}\quad \lVert \mathbf{X}\rVert_{*}\nonumber\\
&&\text{subject to }\quad \mathbf{X}_{\mathbb{O}}=\mathbf{T}_{\mathbb{O}},
\label{opt:nun matrix}
\end{eqnarray}
where the nuclear norm is defined as $\lVert\mathbf{X}\rVert_{*}=\sum_{i=1}^{\operatorname{min}\{M,N\}}\sigma_{i}$, $\sigma_{i}$ is the singular value of $\mathbf{X}$, which can be obtained by singular value decomposition (SVD).

The singular value thresholding (SVT) algorithm can solve  problem (\ref{opt:nun matrix}) by singular value shrinkage operator, which is defined as:
\begin{equation}
\operatorname{SVT}_{\tau}(\mathbf{X})=\mathbf{U}\operatorname{sth}_{\tau}(\Sigma)\mathbf{V}^{\text{T}},\nonumber
\end{equation}
where $\mathbf{U}$ and $\mathbf{V}$ are the right singular vectors and
left singular vectors, respectively. $\operatorname{sth}_{\tau}$ is the well-known soft thresholding operator as follows:
\begin{equation}
\operatorname{sth}_{\tau}=\operatorname{sgn}(x)\text {max}(\lvert x\rvert-\tau,0).
\label{operator: sth}\nonumber
\end{equation}
With $\tau>0$ and a sequence of scalar step sizes ${\delta}\geq1$, $\mathbf{Y}^{0}$, the algorithm defines~\cite{20}:
\begin{eqnarray}
&&\mathbf{X}^{k}=\operatorname{SVT}_{\tau}(\mathbf{Y}^{k-1})\nonumber\\
&&\mathbf{Y}^{k}=\mathbf{Y}^{k-1}+\delta(\mathbf{T}_{\mathbb{O}}-\mathbf{X}_{\mathbb{O}}^{k}),
\end{eqnarray}
until a stopping criterion is reached.

Another method to tackle the problem (\ref{opt:nun matrix}) is that add an additional variable matrix $\mathbf{M}$ and obtain the following equivalent formulation~\cite{27}:
\begin{eqnarray}
&& \minimize_{\mathbf{M}}\quad \lVert \mathbf{M}\rVert_{*}\nonumber\\
&&\text{subject to }\quad \mathbf{M}=\mathbf{X},\quad \mathbf{X}_{\mathbb{O}}=\mathbf{T}_{\mathbb{O}} .
\label{equ: add varia}
\end{eqnarray}
Then we define the following augment Lagrangian function:
\begin{eqnarray}
&&\minimize_{\mathbf{M}}\quad \lVert\mathbf{M}\rVert_{*}+\langle \Lambda,\mathbf{M}-\mathbf{X} \rangle+\frac{\beta}{2}\lVert\mathbf{M}-\mathbf{X}\rVert_\text{F}^{2}\nonumber\\
&& \text{subject to }\quad  \mathbf{X}_{\mathbb{O}}=\mathbf{T}_{\mathbb{O}},
\end{eqnarray}
and according to the alternating direction method of multipliers (ADMM)~\cite{41} framework, one can iteratively update $\mathbf{M}$, $\mathbf{X}$, and $\Lambda$, with the initial values $\mathbf{X}_{\mathbb{O}}^{0}=\mathbf{T}_{\mathbb{O}}$ and $\mathbf{X}_{\bar{\mathbb{O}}}=0$ , the algorithm can be formulated as:
\begin{eqnarray}
&&\mathbf{M}^{k}=\operatorname{D}_{\frac{1}{\beta}}(\mathbf{X}^{k-1}+\frac{1}{\beta}\Lambda^{k-1})\nonumber\\
&& \mathbf{X}_{\bar{\mathbb{O}}}^{k}=(\mathbf{M}^{k}-\frac{1}{\beta}\Lambda^{k-1})_{\bar{\mathbb{O}}}\nonumber\\
&& \Lambda^{k}=\Lambda^{k-1}-\beta(\mathbf{M}^{k}-\mathbf{X}^{k}),
\end{eqnarray}
until a stopping criterion is reached.

\subsection{low rank matrix decomposition based matrix completion}
It's known that SVD becomes increasingly costly as the sizes of the underlying matrices increase, so nuclear norm based algorithms should bear the computational cost required by SVD. It is therefore desirable to exploit an alternative approach that avoids SVD computation all together, by replacing it with
some less expensive computation. Hence, a non-SVD approach in order to more efficiently solve large-scale matrix completion
problems has been proposed.

Based on the low rank matrix decomposition model, the problem (\ref{equ:org}) can be converted into~\cite{24}:
\begin{eqnarray}
&&\minimize_{\mathbf{X},\mathbf{Y},\mathbf{Z}}\quad\frac{1}{2}\lVert\mathbf{X}\mathbf{Y}-\mathbf{Z}\rVert_\text{F}^{2}\nonumber\\
&&\text{subject to} \quad \mathbf{Z}_{\mathbb{O}}=\mathbf{T}_{\mathbb{O}},
\label{equ:XY updated}
\end{eqnarray}
where $\mathbf{X} \in \mathbb{R}^{m\times S}$, $\mathbf{Y} \in \mathbb{R}^{S\times n}$, $\mathbf{Z} \in \mathbb{R}^{m\times n}$, and the integer $S$ is the rank of matrix $\mathbf{Z}$. This problem can be solved by updating these three variables by minimizing ($\ref{equ:XY updated}$) with respect to each one separately while fixing the other two. The updated algorithm is following iterative scheme:
\begin{eqnarray}
&&\mathbf{X}^{k}=\mathbf{Z}^{k-1}(\mathbf{Y}^{k-1})^\text{T}(\mathbf{Y}^{k-1}(\mathbf{Y}^{k-1})^\text{T})^{\dagger}\nonumber\\
&&\mathbf{Y}^{k}=((\mathbf{X}^{k})^\text{T}\mathbf{X}^{k})^{\dagger}((\mathbf{X}^{k})^\text{T}\mathbf{Z}^{k})\nonumber\\
&&\mathbf{Z}^{k}=(\mathbf{X}^{k}\mathbf{Y}^{k})_{\bar{\mathbb{O}}}+\mathbf{T}_{\mathbb{O}}
\end{eqnarray}

\section{Tensor Completion}
\label{sec:4}
Tensor is the generalization of the matrix. Tensor completion is defined as a problem of completing an $N$-th order tensor $ \mathcal{X}\in {\mathbb{R}^{I_1\times I_2\cdots\times I_N}} $ from its known entries given by an index set $\mathbb{O}$. In recent literature, the successful recovery of tensor completion mainly relies on its low rank assumption. The methods for tensor completion include two kinds of approaches. One is based on a given rank and update factors in tensor completion. The one directly minimizes the tensor rank and updates the low-rank tensor.

\begin{figure}[htbp]
	\centering
	\includegraphics[width=300pt, keepaspectratio]{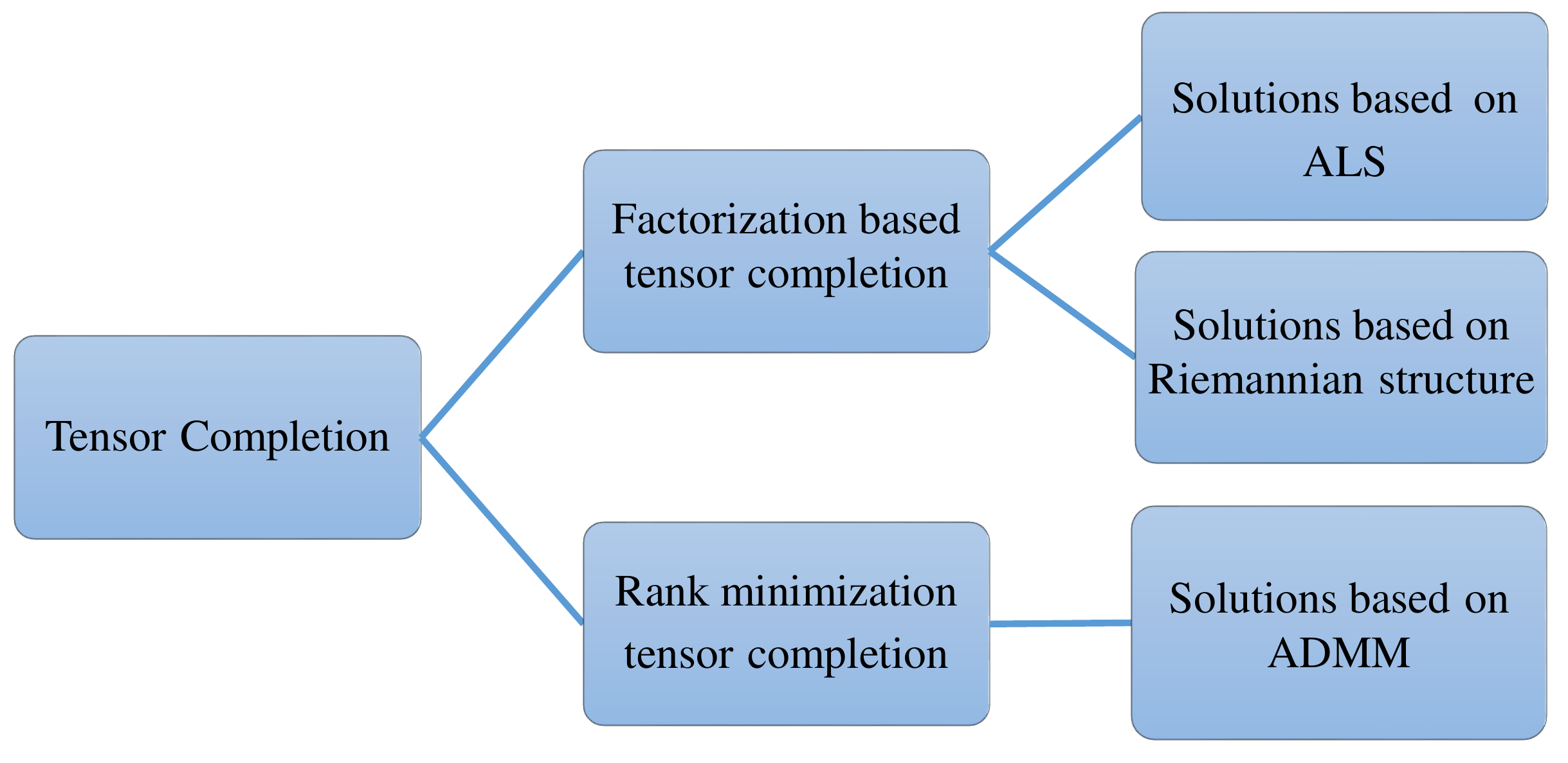}
	\caption{The outline of this chapter. }
	\label{structure}	
\end{figure}

The main framework of this chapter can be illustrated by Fig. \ref{structure}.

\subsection{Factorization based approaches}
In this section, we demonstrate the tensor completion based on a given rank bound with different tensor decomposition methods in detail. The optimization problem of tensor completion with a given rank bound can be formulated as:
\begin{eqnarray}
&&\minimize_{\mathcal{X}} \quad \frac{1}{2}\lVert \mathcal{P}_{\mathbb{O}}(\mathcal{X}-\mathcal{T})\rVert_\text{F}^{2}\nonumber\\
&&\text{subject to} \quad\operatorname{rank(\mathcal{X})}=S,
\label{factorization}
\end{eqnarray}
where $\mathcal{X}$ is the recovered tensor, and $\mathcal{T}$ is the observed tensor. $\operatorname{rank}( \mathcal{X})$  has different forms,  such CP rank, Tucker rank, tensor train rank, tensor ring rank, etc. $S$ is a given bound rank of low rank $ \mathcal{X}$, and $\mathcal{P}_\mathbb{O}$ denotes the random sampling operator, which is defined by:

$$\mathcal{P}_{\mathbb{O}}(\mathcal{A})=\left\{
\begin{array}{rcl}
\mathcal{A}_{i_{1},i_{2},\dots,i_{N}}     &      & (i_{1},i_{2},\dots,i_{N})\in \mathbb{O}  \\
0     &      & \text{otherwise}.
\end{array} \right.$$

\subsubsection{ CP factorization based approaches}
Given a tensor $\mathcal{X} \in \mathbb{R}^{I_{1}\times I_{2}\dots\times I_{N}}$ with a known CP rank, which denotes the number of rank-one tensor, the (\ref{factorization}) can be formulated as the following optimization problem~\cite{42}:
\begin{eqnarray}
&&\minimize_{\mathcal{X}, \mathbf{B}^{(n)}} \quad \frac{1}{2}\lVert \mathcal{P}_{\mathbb{O}}(\mathcal{X}-\mathcal{T})\rVert_\text{F}^{2}\nonumber\\
&&\text{subject to} \quad
\mathcal{X}=\sum_{s=1}^{S} \lambda_{s}\mathbf{b}_{s}^{(1)}\circ\mathbf{b}_{s}^{(2)}\circ\dots\circ\mathbf{b}_{s}^{(N)}.
\end{eqnarray}
This problem can be solved by the alternating least squares (ALS), which updates one variable with other variables fixed~\cite{42}. However, the model is non-convex, which would suffer from the local-minimum and need a good initialization to perform well~\cite{43}. Besides, some other related methods such as CP weighted optimization (CPWOPT)~\cite{44}, Geometric nonlinear conjugate gradient (geomCG)~\cite{45} are proposed to deal with this problem.

Recently, some methods can improve the recovery results through the data prior information. As far as we know, there are mainly two methods that consider the data patterns. One is fully Bayesian CANECOMP/PARAFAC (FBCP) method \cite{46}, which uses the Bayesian inference to find an appropriate tensor rank. The joint distribution of this framework can be written as:
\begin{equation}
p(\mathcal{T}_{\mathbb{O}},\Theta)=p(\mathcal{T}_{\mathbb{O}}|\{\mathbf{B}^{(n)}\}_{n=1}^{N},\tau)\prod_{n=1}^{N}p(\mathbf{B}^{(n)}|\bm{\lambda})p(\bm{\lambda})p(\tau)
\end{equation}
where  the hyperparameters $\Theta=\{B^{(1)}, B^{(2)},\dots, B^{(N)}, \bm{\lambda,\tau}\}$, $\bm{\lambda}=[\lambda_{1}, \lambda_{2},\dots,\lambda_{S}]$, $\tau$ denotes the noise precision. The hyperprior over $\bm{\lambda}$ can be defined by:
\begin{equation}
p(\bm{\lambda})=\prod_{s=1}^{S} \Gamma(\lambda_{s}|a,b),\nonumber
\end{equation}
where $\Gamma(x|a,b)=\frac{b^{a}x^{a-1}e^{-bx}}{\Gamma(a)}$ is a Gamma distribution.

The other is smooth PARAFAC tensor completion (SPC) method \cite{47}, which can be formulated as:
\begin{eqnarray}
&&\minimize_{g_1,\cdots,g_S,\mathbf{u}^{(1)},\ldots,\mathbf{u}^{(N)}} \frac{1}{2}\Arrowvert \mathcal{X}-\mathcal{Y}\Arrowvert_\text{F}^{2}+\sum_{s=1}^{S}\frac{g_{s}^{2}}{2}\sum_{n=1}^{N}\rho^{(n)}\Arrowvert\mathbf{L}^{(n)}\mathbf{u}_{s}^{(N)}\Arrowvert_{p}^{p},\nonumber\\
&&\text{subject to}~ \mathcal{Y}=\sum_{s=1}^{S}g_{s} \mathbf{u}_{s}^{(1)}\circ\mathbf{u}_{s}^{(2)}\circ\cdots\circ\mathbf{u}_{s}^{(n)}, ~\mathcal{X}_{\mathbb{O}}=\mathcal{Z}_{\mathbb{O}},~ \mathcal{X}_{\mathbb{O}}=\mathcal{Y}_{\mathbb{O}},\nonumber\\
&& \Arrowvert \mathbf{u}_{s}^{(N)}\Arrowvert_{2}=1, ~\forall s\in\{1,\ldots,S\},~ \forall n\in\{1, \ldots, N\}.
\end{eqnarray}
where $g_s$ is a coefficient of the $ s $-th rank-one factor tensor in CP decompostion, $\mathbf{u}_{s}^{(n)}$ is $r$-th feature vector, $\circ$ is the outer product of vectors, $\Arrowvert \mathbf{L}^{(n)}\mathbf{u}_{s}^{(n)}\Arrowvert_{p}^{p}$ is a smoothness term, which can be divided into two kinds of total variations. When $p=1$, the term is total variation (TV) and when $p=2$, the constraint term is quadratic variation (QV).
\subsubsection{Tucker factorization based approaches}
Given a tensor $\mathcal{X} \in \mathbb{R}^{I_{1}\times I_{2}\dots \times I_{N}}$ with a known Tucker rank $(S_{1}, S_{2},\dots,S_{N})$, the (\ref{factorization})  can be formulated as the following optimization problem~\cite{48}:
\begin{eqnarray}
&&\minimize_{\mathcal{X}, \mathbf{B}^{(n)}, \mathcal{G}} \quad \frac{1}{2}\lVert \mathcal{P}_{\mathbb{O}}(\mathcal{X}-\mathcal{T})\rVert_\text{F}^{2}\nonumber\\
&&\text{subject to} \quad
\mathcal{X}=\mathcal{G}\times_{1}\mathbf{B}^{(1)}\times_{2}\mathbf{B}^{(2)}\dots\times_{N}\mathbf{B}^{(N)},
\end{eqnarray}
where $\mathcal{G}\in\mathbb{R}^{S_{1}\times S_{2}\times\dots\times S_{N}}$ is the core tensor, and the factor matrix $\mathbf{B}^{(n)}\in \mathbb{R}^{I_{n}\times S_{n}}$.
This problem can be solved by the ALS algorithm~\cite{48} and speeded by the high-order orthogonal iteration (HOOI)~\cite{32},which iteratively calculates one column-wise orthogonal factor matrix corresponding to the dominant singular vectors of $\mathbf{X}^{(n)}$. Besides, when the majority of the tensor entries are unknown, some related methods~\cite{45,49,50} are efficient. They exploit the Riemannian structure on the manifold tensors with known tucker rank and use nonlinear conjugate gradient descent.

In addition, the generalized model of tensor completion with a predefined Tucker rank is popular and can be formulated as:
\begin{eqnarray}
&&\minimize_{\mathcal{X}} \quad \frac{1}{2}\lVert \mathbf{t}-\mathcal{P}(\mathcal{X})\rVert_\text{F}^{2}\nonumber\\
&& \text{subject to} \quad \operatorname{rank}_{\text{Tucker}}(\mathcal{X})\leq (S_{1}, S_{2},\dots,S_{N}),
\end{eqnarray}
where $\mathbf{t}\in \mathbb{R}^{M}$ is the observation data, $\mathcal{P}: \mathbb{R}^{I_{1}\times I_{2},\dots\times I_{N}} \longmapsto\mathbb{R}^{M}$ is the linear measurement operator.
This problem can be solved by tensor iterative hard thresholding (TIHT)~\cite{51}, which employs the truncated higher order singular value decomposition (HOSVD)~\cite{33} as the thresholding operator and achieves a quasi-optimal low rank tensor approximation. In ~\cite{52}, an improved step size TIHI (ISS-TIHI) is proposed to increase the convergence speed. Besides, other related works such as sequentially
optimal modal projections (SeMP)~\cite{53}, sequential rank-one approximation
and projection (SeROAP)~\cite{54}, and sequential low-rank approximation and projection (SeLRAP)~\cite{55} are proposed to effectively resolve tensor recovery problems.
\subsubsection{Tensor train factorization based approaches}
Given a tensor $\mathcal{X} \in \mathbb{R}^{I_{1}\times I_{2}\dots\times I_{N}}$ with a known tensor train rank $(S_{1}, S_{2},\dots,S_{N}, S_{N+1})$, the (\ref{factorization}) can be formulated as the following optimization problem~\cite{56}:
\begin{eqnarray}
&&\minimize_{\mathcal{X},  \mathcal{U}} \quad \frac{1}{2}\lVert\mathcal{X}-\mathcal{Y} \rVert_\text{F}^{2}\nonumber\\
&&\text{subject to} \quad
\mathcal{P}_{\mathbb{O}}(\mathcal{X})=\mathcal{P}_{\mathbb{O}}(\mathcal{T}), \quad \mathcal{Y} \in \mathbb{T}_\text{TT}(S_{1}, S_{2},\dots,S_{N}).
\end{eqnarray}
where $\mathcal{Y} \in \mathbb{T}_\text{TT}(S_{1}, S_{2},\dots,S_{N}) $ means $
\mathcal{Y}(i_1,i_2,\cdots, i_N)=\mathcal{U}_{1}(:,i_1,:)\mathcal{U}_{2}(:,i_2,:)\cdots \mathcal{U}_{N}(:,i_N,:)$, $\mathcal{U}_{n}\in\mathbb{R}^{S_{n}\times I_{n}\times S_{n+1}}$. This problem is non-convex and exists a quasi-optimal approximation solution. It can be calculated using nonlinear block Gauss-Seidel iteration with respect to the blocks $\mathcal{U}_{1}, \mathcal{U}_{2},\dots, \mathcal{U}_{N},\mathcal{X}$~\cite{56, 57}. For each small block, it can be solved by least squares fit and this approach can accelerate the convergence by successive over-relaxation, which is called as alternating directions fitting ~\cite{56}. Besides, with fixed tensor train rank, this problem can be solved using a nonlinear conjugate gradient scheme within the framework of Riemannian optimization\cite{58}, which can make the storage complexity scaling linear with the number of dimensions.

In addition, with the success of matrix completion~\cite{59} and tensor completion using Tucker rank~\cite{60, 61}, a similarly way of tensor completion using fixed tensor train rank bound was proposed in the form of~\cite{62}:
\begin{eqnarray}
&& \minimize_{\mathbf{U}_{n},\mathbf{V}_{n},\mathcal{X}}\quad \sum_{n=1}^{N-1}\frac{\alpha_{n}}{2}\lVert\mathbf{U}_{n}\mathbf{V}_{n}-\mathbf{X}_{<n>}\rVert_\text{F}^{2}\nonumber\\
&& \text{subject to}\quad \mathcal{P}_{\mathbb{O}}(\mathcal{X})=\mathcal{P}_{\mathbb{O}}(\mathcal{T})
\end{eqnarray}
where $\mathbf{X}_{<n>}\in\mathbb{R}^{\prod_{j=1}^{n}I_{j}\times \prod_{j=n+1}^{N}I_{j}}$. With a given rank $S_{n}$, the matrix $\mathbf{X}_{<n>}$ can be factorized as $\mathbf{X}_{<n>}=\mathbf{U}_{n}\mathbf{V}_{n}$ with $\mathbf{U}_{n}\in\mathbb{R}^{\prod_{j=1}^{n}I_{j}\times S_{n}}$ and $\mathbf{V}_{n}\in\mathbb{R}^{S_{n}\times \prod_{j=n+1}^{N}I_{j}}$. For this problem, the block coordinate descent (BCD) method~\cite{63} was used to optimize to alternatively optimize different groups of variables ($\mathbf{U}_{n},\mathbf{V}_{n},\mathcal{X}$).
\subsubsection{Tensor Ring factorization based approaches}
To the best of our knowledge, low rank tensor ring completion~\cite{64} is the only one work for tensor completion with fixed tensor ring rank $(S_{1}, S_{2},\dots,S_{N})$. Given a tensor  $\mathcal{X} \in \mathbb{R}^{I_{1}\times I_{2}\dots\times I_{N}}$ under tensor ring model, the optimization problem can be formulated as:
\begin{eqnarray}
&&\minimize_{\mathcal{U}} \quad \frac{1}{2}\lVert\mathcal{P}_{\mathbb{O}}(\mathcal{X}-\mathcal{T}) \rVert_\text{F}^{2}\nonumber\\
&&\text{subject to}
\quad \mathcal{X} \in \mathbb{T}_\text{TR}(S_{1}, S_{2},\dots,S_{N}).
\label{TR_ALS}
\end{eqnarray}
where $\mathcal{X} \in \mathbb{T}_\text{TR}(S_{1}, S_{2},\dots,S_{N}) $ means $
\mathcal{X}(i_1,i_2,\cdots, i_N)=\operatorname{tr}(\mathcal{U}_{1}(:,i_1,:)\mathcal{U}_{2}(:,i_2,:)\cdots \mathcal{U}_{N}(:,i_N,:))$, $\mathcal{U}_{n}\in\mathbb{R}^{S_{n}\times I_{n}\times S_{n+1}}$.
The tensor train ranks present the distribution of the large for the middle factors and the lower for the boarder factors. And in fact, the closer to the boarder factor, the lower the tensor train ranks are. The tensor train decomposition drawbacks can be alleviate by tensor ring factorization. The optimization problem (\ref{TR_ALS}) can be solved by ALS with a good initial point choosing. In~\cite{64}, the author proposed a novelty initialization algorithm which called tensor ring approximation to make sure the ALS algorithm performs well.

\subsubsection{t-SVD factorization based approaches}
In this section, we consider a tensor $\mathcal{X}\in\mathbb{R}^{I_{1}\times I_{2}\times I_{3}}$ with a fixed tubal rank. The optimization problem of low tubal rank tensor completion problem can be formulated as:
\begin{eqnarray}
&&\minimize_{\mathcal{X}\in\mathbb{R}^{I_{1}\times I_{2}\times I_{3}}} \quad \lVert \mathcal{P}_{\mathbb{O}}(\mathcal{X}-\mathcal{T})\rVert_\text{F}\nonumber\\
&&\text{subject to} \quad\operatorname{rank(\mathcal{X})}\leq S,
\label{tsvd_B}
\end{eqnarray}
This problem can be solved by decomposing the target tensor as the circular convolution of two low tubal rank tensor, which can be rewritten as~\cite{65}:
\begin{eqnarray}
&&\minimize_{\mathcal{X}, \mathcal{A}, \mathcal{B}} \quad \lVert \mathcal{P}_{\mathbb{O}}(\mathcal{A}\ast\mathcal{B}^{\dagger}-\mathcal{T})\rVert_\text{F}\nonumber\\
&&\text{subject to} \quad\mathcal{X}=\mathcal{A}\ast\mathcal{B}^{\dagger}
\end{eqnarray}
where $\mathcal{A}\in\mathbb{R}^{I_{1}\times S\times I_{3}}, \mathcal{B}\in\mathbb{R}^{I_{2}\times S\times I_{3}}$. This low tubal rank tensor completion model can be solved by alternating minimization algorithm \cite{24, 66}.

Furthermore, a tensor-CUR decomposition based model was proposed to solve the problem (\ref{tsvd_B}), which can be converted to the following optimization problem~\cite{67}:
\begin{eqnarray}
&&\minimize_{\mathcal{U}} \quad \lVert \mathcal{P}_{\mathbb{O}}(\mathcal{C}\ast\mathcal{U}\ast\mathcal{R}-\mathcal{T})\rVert_\text{F}\nonumber\\
&&\text{subject to} \quad\mathcal{X}=\mathcal{C}\ast\mathcal{U}\ast\mathcal{R},\quad \mathcal{P}_{\mathbb{O}}(\mathcal{X})=\mathcal{P}_{\mathbb{O}}(\mathcal{T})
\end{eqnarray}
where $\mathcal{C}\in\mathbb{R}^{I_{1}\times I_{1}\times I_{3}}$, $\mathcal{U}\in \mathbb{R}^{I_{1}\times I_{2}\times I_{3}}$, $\mathcal{R}\in \mathbb{R}^{I_{2}\times I_{2}\times I_{3}}$. It can be solved by the matrix-CUR approximation for each frontal face. Besides, In \cite{68}, the authors show the smooth manifold of a fixed tubal rank tensor can be seen as an element of the product manifold of fixed low-rank matrices in
the Fourier domain. This problem can be solved by the traditional convex methods such as conjugate gradient descent~\cite{69}.
\subsection{Rank minimization model}
In practice, the tensor rank bounds may not be available in some applications. When only a few observations are available, the choice of high rank bound may led to over-fitting. To avoid the occurrence of it, another group of methods is to directly minimize the tensor rank, which can be given as follows:
\begin{eqnarray}
&&\minimize_{\mathcal{X}}\quad\operatorname{rank}(\mathcal{X})\nonumber\\
&&\text{subject to}\quad \mathcal{X}_{\mathbb{O}}=\mathcal{T}_{\mathbb{O}},
\label{rank_minimize}
\end{eqnarray}
where $\mathcal{X}_{\mathbb{O}}=\mathcal{T}_{\mathbb{O}}$ means the $\mathcal{X}_{i_{1},\dots,i_{N}}=\mathcal{T}_{i_{1},\dots,i_{N}}$, the index $\{i_{1},\dots,i_{N}\}$ is in the observation index $\mathbb{O}$, where $\operatorname{rank}( \mathcal{X})$ denotes the rank of tensor variable $ \mathcal{X} $. There are different kinds of tensor ranks, such as CP rank, Tucker rank, tubal rank, tensor train rank, etc. With different definitions of tensor rank, there are many methods optimization models for tensor completion problems. However, the rank is a non-convex function with respect to $\mathcal{X}$, and the problem (\ref{rank_minimize}) is NP-hard~\cite{70}. Most existing methods~\cite{71, 72} are using nuclear norm (trace norm) as the convex surrogate of non-convex rank function. In particular, the tensor nuclear norm for CP rank is impossible to solve and  to the best of our knowledge, tensor nuclear norm for tensor ring rank has never been investigated for tensor completion. Hence, we will mainly introduce the methods of Tucker rank minimization, tubal rank minimization and tensor train rank minimization for tensor completion in the following subsections.

\subsubsection{Tucker rank minimization model}
Given a tensor $\mathcal{X}\in\mathbb{R}^{I_{1}\times I_{2}\dots\times I_{N}}$, the tensor completion based on minimizing Tucker rank can be formulated as:
\begin{eqnarray}
&&\minimize_{\mathcal{X}}\quad\operatorname{rank}_{\text{Tucker}}(\mathcal{X})\nonumber\\
&&\text{subject to}\quad \mathcal{X}_{\mathbb{O}}=\mathcal{T}_{\mathbb{O}},
\label{rankT_minimize}
\end{eqnarray}
where $\operatorname{rank}_{\text{Tucker}}(\mathcal{X})=(\operatorname{rank}(\mathbf{X}_{[1]}), \operatorname{rank}(\mathbf{X}_{[2]}),\dots,\operatorname{rank}(\mathbf{X}_{[N]})$~\cite{32}, and $\operatorname{rank}(\mathbf{X}_{[n]})$ denotes the rank of the unfolding matrix $\mathbf{X}_{[n]}$. Under the definition of Tucker rank, the optimization problem (\ref{rankT_minimize}) can be written as~\cite{73, 74, 60}:
\begin{eqnarray}
\label{eq: Tucker rank minimization}
&&\minimize_{\mathbf{X}_{[n]}}\quad \sum^N_{n=1} w_n \operatorname{rank}(\mathbf{X}_{[n]}) \nonumber\\
&& \text {subject to} \quad \mathcal X_{\mathbb{O}}= \mathcal T _{\mathbb{O}},
\end{eqnarray}
where $ w_n, n=1, \cdots, N $ are the weights with  $ \sum^N_{n=1} w_n=1 $. In fact, the matrix rank functions in the optimization model (\ref{eq: Tucker rank minimization}) is nonconvex, but it can be relaxed to the matrix nuclear norm approximation. This generates the following optimization problem~\cite{73}:
\begin{eqnarray}
\label{eq: Nuclear norm minimization}
&&\minimize_{\mathbf{X}_{[n]}}\quad \sum^N_{n=1} w_n\Arrowvert \mathcal{X}_{[n]}\Arrowvert_*\nonumber\\
&&\text {subject to} \quad \mathcal{X}_{\mathbb{O}}= \mathcal{T}_{\mathbb{O}}.	
\end{eqnarray}
Where $\Arrowvert \mathcal{X}_{[n]}\Arrowvert_* $ is the tightest convex envelop for $\operatorname{rank}(\mathbf{X}_{[n]})$. The problem (\ref{eq: Nuclear norm minimization}) can be solved using simple low-rank tensor completion (SiLRTC) and strictly solved using high accuracy low-rank tensor completion (HaLRTC) by adding an equation constraint~\cite{73}. In order to improve the performance, the method which uses volume measurement~\cite{75} is proposed and can be formulated as:
\begin{eqnarray}
&& \minimize_{\mathcal{X}} \quad\lVert \mathcal{X}\rVert_{\text{vol}}\nonumber\\
&&\text{subject to} \quad\mathcal{X}_{\mathbb{O}}=\mathcal{T}_{\mathbb{O}},
\end{eqnarray}
where the volume of tensor is defined as:
\begin{eqnarray}
\lVert\mathcal{X}\rVert_{\text{vol}}&=&\prod_{i=1}^{N}\lVert\mathbf{X}_{[i]}\rVert_{\text{vol}}\nonumber\\
&=&\prod_{i=1}^{N}\prod_{j=1}^{N_{i}}\sigma_{j}(\mathbf{X}_{[i]}),\nonumber
\end{eqnarray}
$\sigma_{j}(\mathbf{X}_{[i]})$ is the $j$-th singular value of mode-i unfolding matrix. This model can be solved by the ADMM. However, the sum of weighted nuclear norms model may be sub-optimal with the dimension increasing~\cite{76, 77}. To handle the issue, a more appropriate convex model, which makes the mode-n unfolding matrix more balanced and maintains the low rank property was proposed~\cite{76}. In addition, the core tensor
trace-norm minimization (CTNM)~\cite{78} method using ADMM algorithm can alleviate computational cost in large scale problems.

Furthermore, in order to enhance the recovery quality for some real data, such as natural images, videos, magnetic resonance imaging (MRI) images and hyperspectral images, investigating the intrinsic structure of data is necessary. Hence, some smooth priors can be added to improve the recovery quality for images completion~\cite{79,80,81}. There are two groups of state-of-the-art models for tensor completion considering the prior information.

Based on the factor prior, one presents the nuclear norm of individual factor matrices as the following optimization problem:
\begin{eqnarray}
\label{STDC}
&&\minimize_{\mathcal{X},\mathcal{G},\mathbf{U}^{(n)}} \quad \sum_{n=1}^{N}w_{n}\lVert\mathbf{U}^{(n)}\rVert_{*}+\beta \operatorname{tr}((\mathbf{U}_{1}\otimes\dots\otimes\mathbf{U}_{N})\mathbf{L}(\mathbf{U}_{1}\otimes\dots\otimes\mathbf{U}_{N})^{\text{T}})+\gamma\lVert\mathcal{G}\rVert_\text{F}^{2}\nonumber\\
&&\text{subject to}\quad \mathcal{X}=\mathcal{G}\times_{1}\mathbf{U}_{1}^{\text{T}}\dots\times_{n}\mathbf{U}_{n}^{\text{T}},\quad \mathcal{X}_{\mathbb{O}}=\mathcal{T}_{\mathbb{O}}.
\end{eqnarray}
where $\beta \operatorname{tr}((\mathbf{U}_{1}\otimes\dots\otimes\mathbf{U}_{N})\mathbf{L}(\mathbf{U}_{1}\otimes\dots\otimes\mathbf{U}_{N})^{\text{T}})+\gamma\lVert\mathcal{G}\rVert_\text{F}^{2}$ are the regularization items, $\beta$ and $\gamma$ are the weighting parameters between nuclear norm and regularization items. $\mathbf{L}$ is a matrix designed by the prior information, and $ \operatorname{tr}((\mathbf{U}_{1}\otimes\dots\otimes\mathbf{U}_{N})\mathbf{L}(\mathbf{U}_{1}\otimes\dots\otimes\mathbf{U}_{N})^{\text{T}})$ can be interpreted as the constraint for the local similarity of visual data. This problem can be solved using augmented Lagrange multiplier (ALM) algorithm \cite{82}.

Based on the recovered tensor prior, the other proposes the simultaneous minimization of the tensor nuclear norm and total variation (TV). The optimization model is as follow:
\begin{eqnarray}
\label{LRTV_PDS}
&&\minimize_{\mathcal{X}} \quad\lVert (1-\rho)\Arrowvert\mathcal{X}\Arrowvert_{*}+\rho\Arrowvert\mathcal{X}\Arrowvert_{\text{TV}}\nonumber\\
&&\text{subject~to}~\nu_\text{min}\leq\mathcal{X}\leq\nu_\text{max},~\Arrowvert\mathcal{Z}_{\mathbb{O}}-\mathcal{X}_{\mathbb{O}}\Arrowvert_\text{F}^{2}\leq\delta.
\end{eqnarray}
where $\rho$ is the weight between nuclear norm and TV terms, and $ \delta $ is the error bound.The first inequality constraint in (\ref{LRTV_PDS}) is to impose all values of output tensor within a given range [$\nu_\text{min}$,$\nu_\text{max}$]. The second inequality constraint implies that the recovered data entries and the observed data entries in the set $\mathbb{O}$ are generally consistent but allows existing a certain amount of noise. The problem of simultaneous minimization of low rank and TV terms can be solved by the primal-dual splitting method (LRTV-PDS)~\cite{80,81}.

\subsubsection{Tensor train rank minimization model}
Given a tensor $\mathcal{X}\in\mathbb{R}^{I_{1}\times I_{2}\dots\times I_{N}}$, the tensor completion based on minimizing tensor train rank can be formulated as:
\begin{eqnarray}
&&\minimize_{\mathcal{X}}\quad\operatorname{rank}_{\text{TT}}(\mathcal{X})\nonumber\\
&&\text{subject to}\quad \mathcal{X}_{\mathbb{O}}=\mathcal{T}_{\mathbb{O}},
\label{rankTT_minimize}
\end{eqnarray}
where $\operatorname{rank}_{\text{TT}}(\mathcal{X})=(\operatorname{rank}(\mathbf{X}_{<1>}), \operatorname{rank}(\mathbf{X}_{<2>}),\dots,\operatorname{rank}(\mathbf{X}_{<N-1>})$~\cite{36}, and $\operatorname{rank}(\mathbf{X}_{<n>})$ denotes the rank of the unfolding matrix $\mathbf{X}_{<n>}$, so it can be written as (\ref{rankTT_minimize}):
\begin{eqnarray}
\label{eq: Tensor train rank minimization}
&&\minimize_{\mathbf{X}_{<n>}}\quad \sum^{N-1}_{n=1} w_n \operatorname{rank}(\mathbf{X}_{<n>}) \nonumber\\
&& \text {subject to} \quad \mathcal X_{\mathbb{O}}= \mathcal T _{\mathbb{O}},
\end{eqnarray}
where $ w_n$, $n=1, \cdots, N-1 $ are the weights with  $ \sum_{n=1}^{N-1} w_n=1 $.
This optimization problem can be approximated by the matrix nuclear norm as follows~\cite{62}:
\begin{eqnarray}
\label{eq: TT Nuclear norm minimization}
&&\minimize_{\mathbf{X}_{<n>}}\quad \sum^{N-1}_{n=1} w_n\Arrowvert \mathbf{X}_{<n>}\Arrowvert_*\nonumber\\
&&\text {subject to} \quad \mathcal{X}_{\mathbb{O}}= \mathcal{T}_{\mathbb{O}}.	
\end{eqnarray}

Especially, the unfolding $\mathbf{X}_{<n>}$ is more balanced than the unfolding $\mathbf{X}_{[n]}$. For this model when $n=\frac{N}{2}$, it is similar to the case in~\cite{76}. The problem (\ref{eq: TT Nuclear norm minimization}) can be solved by simple low-rank tensor completion via tensor train (SiLRTC-TT)~\cite{62}.

\subsubsection{Tubal rank minimization model}
Given a tensor $\mathcal{X}\in\mathbb{R}^{I_{1}\times I_{2}\times I_{3}}$, the tensor completion based on minimizing tensor tubal rank can be formulated as:
\begin{eqnarray}
&&\minimize_{\mathcal{X}}\quad\operatorname{rank}_{\text{tubal}}(\mathcal{X})\nonumber\\
&&\text{subject to}\quad \mathcal{X}_{\mathbb{O}}=\mathcal{T}_{\mathbb{O}},
\label{ranktubal_minimize}
\end{eqnarray}
where the tensor tubal rank, denoted by $\mathrm{rank_{tubal}}(\mathcal{X})$, is the number of nonzero singular tubes of $\mathcal{S}$, where $\mathcal{S}$ is from  $\mathcal{X}=\mathcal{U}* \mathcal{S}*\mathcal{V^\mathrm{T}}$~\cite{83}. Following the optimization approach in~\cite{73}, the problem (\ref{ranktubal_minimize}) is firstly formulated as~\cite{83}:
\begin{eqnarray}
&&\minimize_{\mathcal{X}}\quad \lVert\mathcal{X}\rVert_\text{TNN}\nonumber\\
&&\text{subject to} \quad\mathcal{X}_{\mathbb{O}}=\mathcal{T}_{\mathbb{O}}
\end{eqnarray}
where $\lVert\mathcal{X}\rVert_\text{TNN}$ is the nuclear norm of $\mathcal{X}$, which is equal to the tensor nuclear norm (TNN ) of $\operatorname{blkdiag(\mathcal{\hat{X}})}$.
Here $\text{blkdiag}(\mathcal{\hat{X}})$ is a block diagonal matrix defined as follows:
\begin{equation}
\mbox{blkdiag}(\hat{ \mathcal{X} })=
\left[
\begin{array}{cccc}
\hat{ \mathcal{X} }^{(1)} & & &\\
& \hat{ \mathcal{X} }^{(2)}& &\\
& & \ddots& \\
& & &\hat{ \mathcal{X} }^{(\emph{n}_3)}
\end{array}
\right],\nonumber
\end{equation}
where $\hat{ \mathcal{X} }^{(i)}$ is the $i$-th frontal slice of $\hat{ \mathcal{X} }, i=1, 2,..., I_3$.
Besides, some related works were proposed to complete the video, which has the redundancy in frames and spatial resolution, such as a twist tensor nuclear norm~\cite{84} and TNN~\cite{85}.

Furthermore, to better estimate the tensor tubal rank, some new convex envelopes, such as the tensor truncated nuclear norm (T-TNN) and weighted tensor nuclear norm (W-TNN), were proposed to replace the traditional TNN in t-SVD.

The optimization model based on T-TNN can be rewritten as~\cite{86}:
\begin{eqnarray}
&&\minimize_{\mathcal{X}}\quad \lVert \mathcal{X}\rVert_{*}-\max_{\mathcal{A}\ast\mathcal{A}^\text{T}=\mathcal{I},~ \mathcal{B}\ast\mathcal{B}^\text{T}=\mathcal{I}} \operatorname{tr}(\mathcal{A}\ast\mathcal{X}\ast\mathcal{B}^\text{T})\nonumber\\
&&\text{subject to}\quad \mathcal{X}_{\mathbb{O}}=\mathcal{T}_{\mathbb{O}},
\end{eqnarray}
where $\mathcal{A}\in\mathbb{R}^{I_{1}\times S\times I_{3}}$ and $\mathcal{B}\in\mathbb{R}^{I_{2}\times S\times I_{3}}$. This problem can be divided  into two sub-problems. One updates the variable tensor $\mathcal{X}$ by ADMM with the tensors $\mathcal{A}$ and $\mathcal{B}$ fixed. And the other one updates the variable tensors $\mathcal{A}$ and $\mathcal{B}$ using the updated tensor $\mathcal{X}$, which can be solved by t-SVD in Definition 16.

The another optimization model based on W-TNN is as follow\cite{87}:
\begin{eqnarray}
&&\minimize_{\mathcal{X}}\quad \lVert\mathcal{X},\mathcal{W}_{\mathcal{X}}\rVert_{\mathcal{W}\circledast}\nonumber\\
&&\text{subject to}\quad \mathcal{X}_{\mathbb{O}}=\mathcal{T}_{\mathbb{O}},
\end{eqnarray}
where $\mathcal{W}_{\mathcal{X}}$ is the weight tensor of $\mathcal{X}$, and $\lVert\mathcal{X},\mathcal{W}_{\mathcal{X}}\rVert_{\mathcal{W}\circledast}$ is the weighted tensor nuclear norm operator, which is defined as:
where $\Sigma_{f}(i,i,j)$ is the singular value of $\mathcal{X}$. Through introducing an auxiliary tensor and converting it to Inexact Augmented Lagrangian function, this problem can be solved using Inexact Augmented Lagrange Multiplier (IALM).

\subsection{other variants}
\begin{figure}[htbp]
	\centering
	\includegraphics[width=350pt, keepaspectratio]{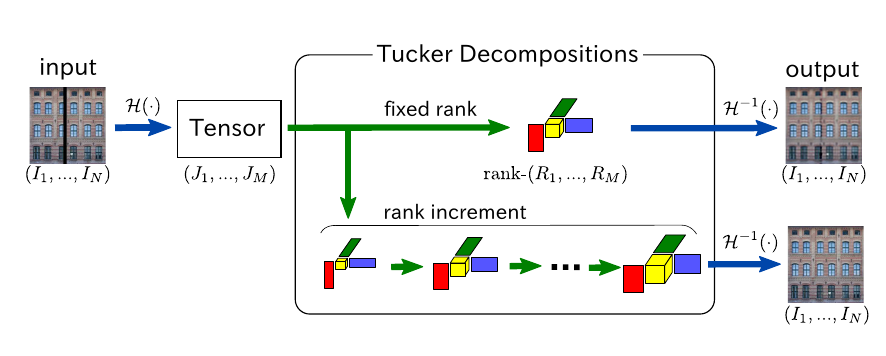}
	\caption{Illustration of the tensor completion using embedding transform.}
	\label{HTuckercompletion}
\end{figure}

Except for the methods that complete the tensor in its original space, extending it to a higher-order tensor is also a good method in some special cases. In~\cite{88}, the authors proposed a method which using some embedding transform to a higher space for the tensor completion. This method suits for the case that all entries are missed in some continuous slices, which is illustrated in Fig. \ref{HTuckercompletion}.

\section{Experiments}
\label{sec:5}
In this section, we conduct several experiments on color image $``$ peppers $"$ using different tensor completion methods, which include factorization based approaches and rank minimization based approaches. The sampling methods of all test experiments are random sampling with respect to different sampling ratios (SR), which is defined as:
\begin{equation}
\text{SR}=\frac{O}{\prod_{n=1}^{N}I_{n}}
\end{equation}
where $O$ is the number of observation entries, $I_{n}$, $n=1,\dots,N$ is the dimension corresponding to mode $k$.

To measure the recovered performance of tensor completion algorithms, the relative error between the original tensor and recovered tensor, which is the most used evaluation metric, has been defined as:
\begin{equation}
\text{Rel}=\lVert \mathcal{\bar{X}}-\mathcal{X}\rVert_\text{F}/\lVert\mathcal{X}\rVert_\text{F},
\end{equation}
where $\mathcal{\bar{X}}$ is the recovered tensor and $\mathcal{X}$ is the original tensor. Besides, peak signal-to-noise ratio (PSNR) and structural similarity index (SSIM) are the commonly applied in images recovery. The PSNR is the error between the original and reconstructed images, which is defined as:
\begin{equation}
\text{PSNR}=20\log_{10}(\operatorname{MAX}/\operatorname {\sqrt{MSE}}) ,
\end{equation}
where $\operatorname{MAX} $ is the maximum value of all the pixels in the image, and mean squared error (MSE) is defined as:
\begin{equation}
\operatorname{MSE}=\lVert\mathcal{P}_{\mathbb{O}}\ast( \hat{\mathcal{X}}-\mathcal{X})\rVert_\text{F}^{2}/\lVert\mathcal{P}_{\mathbb{O}}\rVert_\text{F}^{2}.
\end{equation}
The SSIM measures the similarity of two images based on three comparison measurements with respect to luminance, contrast and structure, which is defined as:
\begin{equation}
\operatorname{SSIM}=l(\mathcal{X},\mathcal{\hat{X}})^{\alpha}\cdot c(\mathcal{X},\mathcal{\hat{X}})^{\beta}\cdot s(\mathcal{X},\mathcal{\hat{X}})^{\gamma},
\end{equation}
where $\alpha$, $\beta$ and $\gamma$ are the weights corresponding to the luminance measurement $l(\mathcal{X},\mathcal{\hat{X}})$, contrast measurement $c(\mathcal{X},\mathcal{\hat{X}})$, structure measurement $s(\mathcal{X},\mathcal{\hat{X}})$.

\begin{figure}[htbp]
	\centering
	\includegraphics[width=300pt, keepaspectratio]{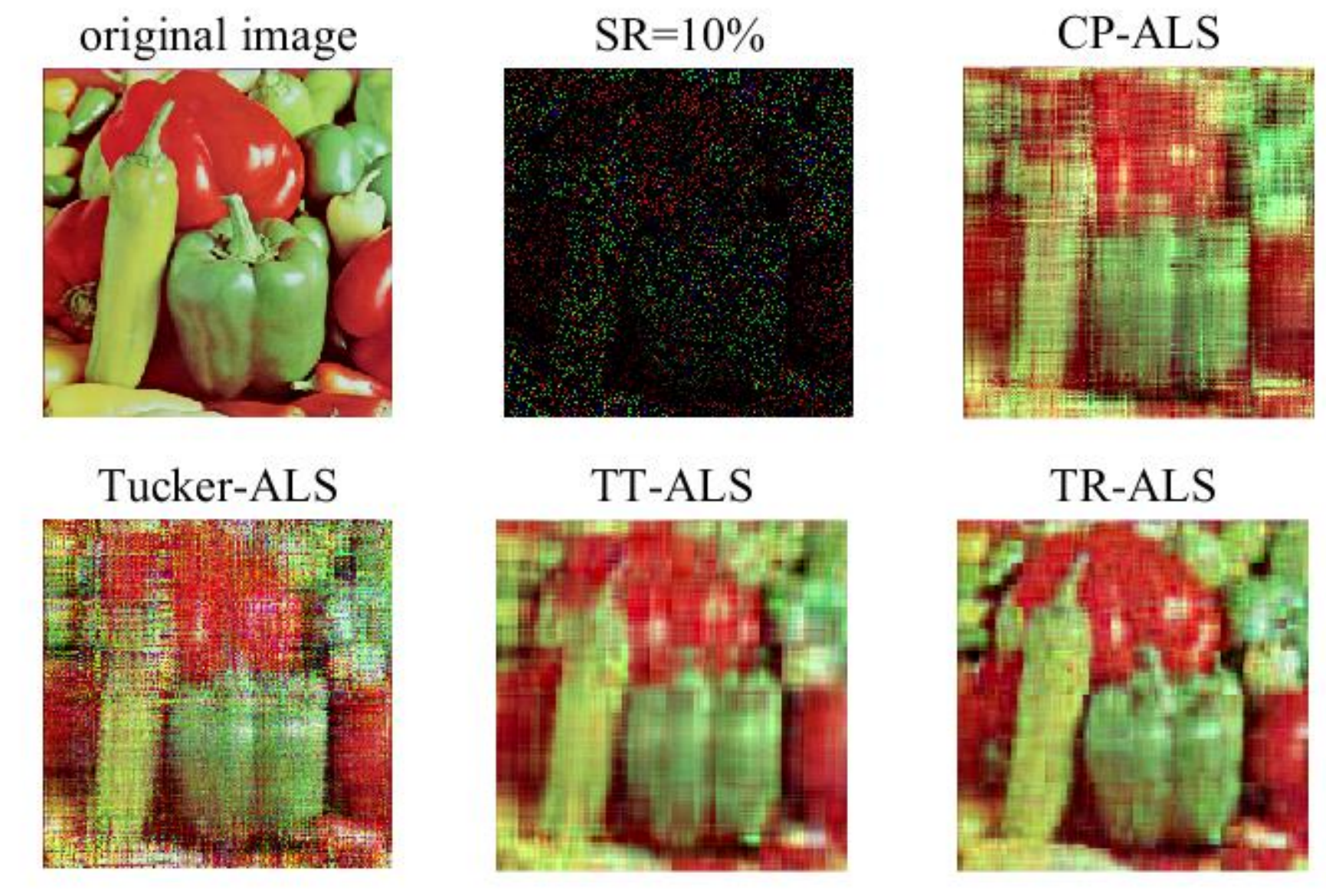}
	\caption{Recovered the peppers image using different methods with a given bound rank.}
	\label{ALS}
\end{figure}

\begin{figure}[htbp]
	\centering
	\includegraphics[width=300pt, keepaspectratio]{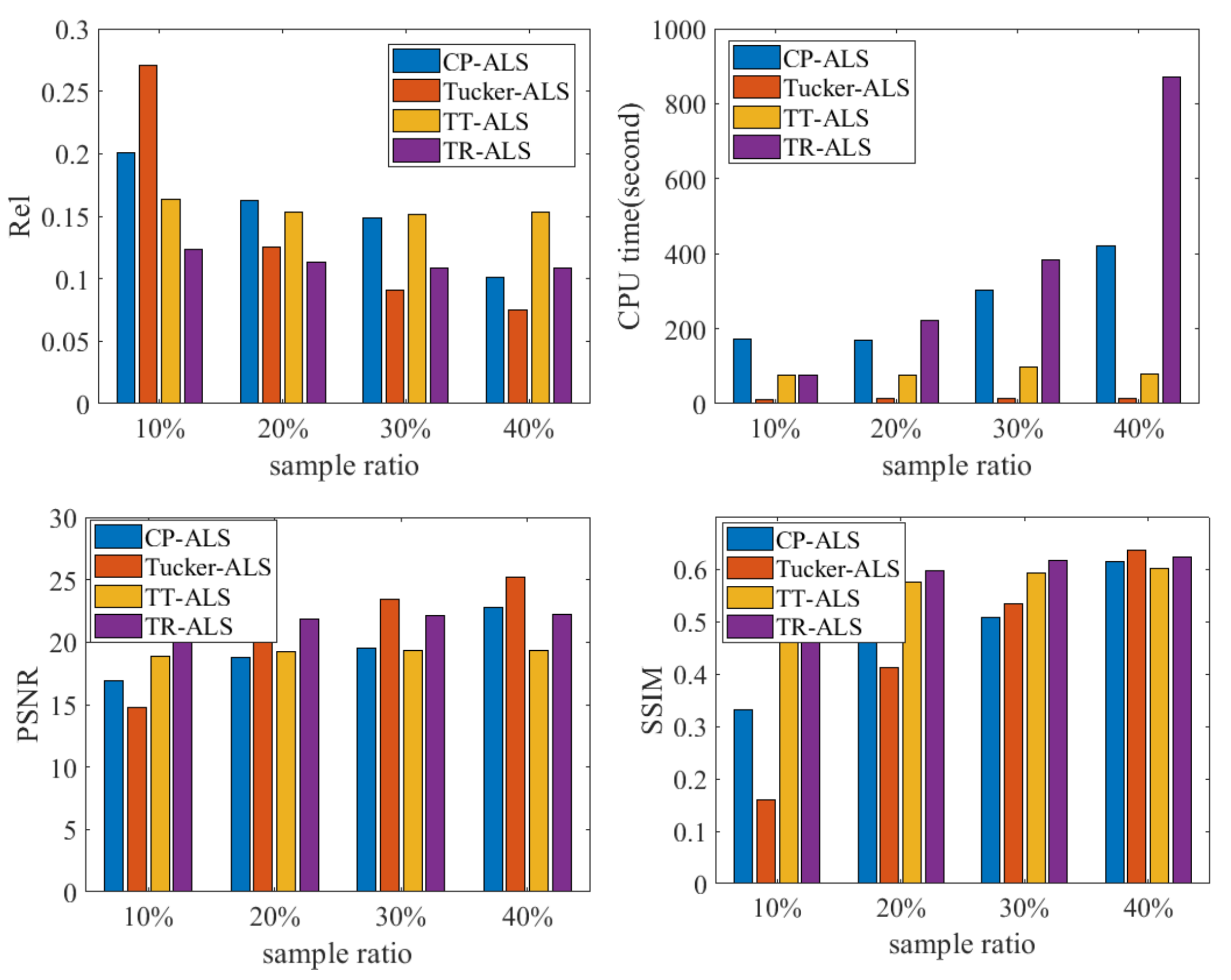}
	\caption{The performance of different decomposition method with a bound rank in terms of Rel, CPU time, PSNR and SSIM.}
	\label{RelALS}
\end{figure}

Fig. \ref{ALS} shows the recovered images using ALS algorithms using different decomposition with a given bound rank. We set the bound ranks parameters as: $10$, $[10,10,3]$, $[1,10,10,1]$, $[10,10,10]$ for CP rank bound, tucker rank bound, tensor train rank bound and tensor ring rank bound, respectively.
Fig. \ref{RelALS} shows the metric using Rel, PSNR and SSIM with SR from 10\% to 40\%. From the Fig. \ref{RelALS}, we could observe that different decomposition methods with the same rank bound under the SR changes from 10\% to 40\%, the result is different. 

\begin{figure}[htbp]
	\centering
	\includegraphics[width=300pt, keepaspectratio]{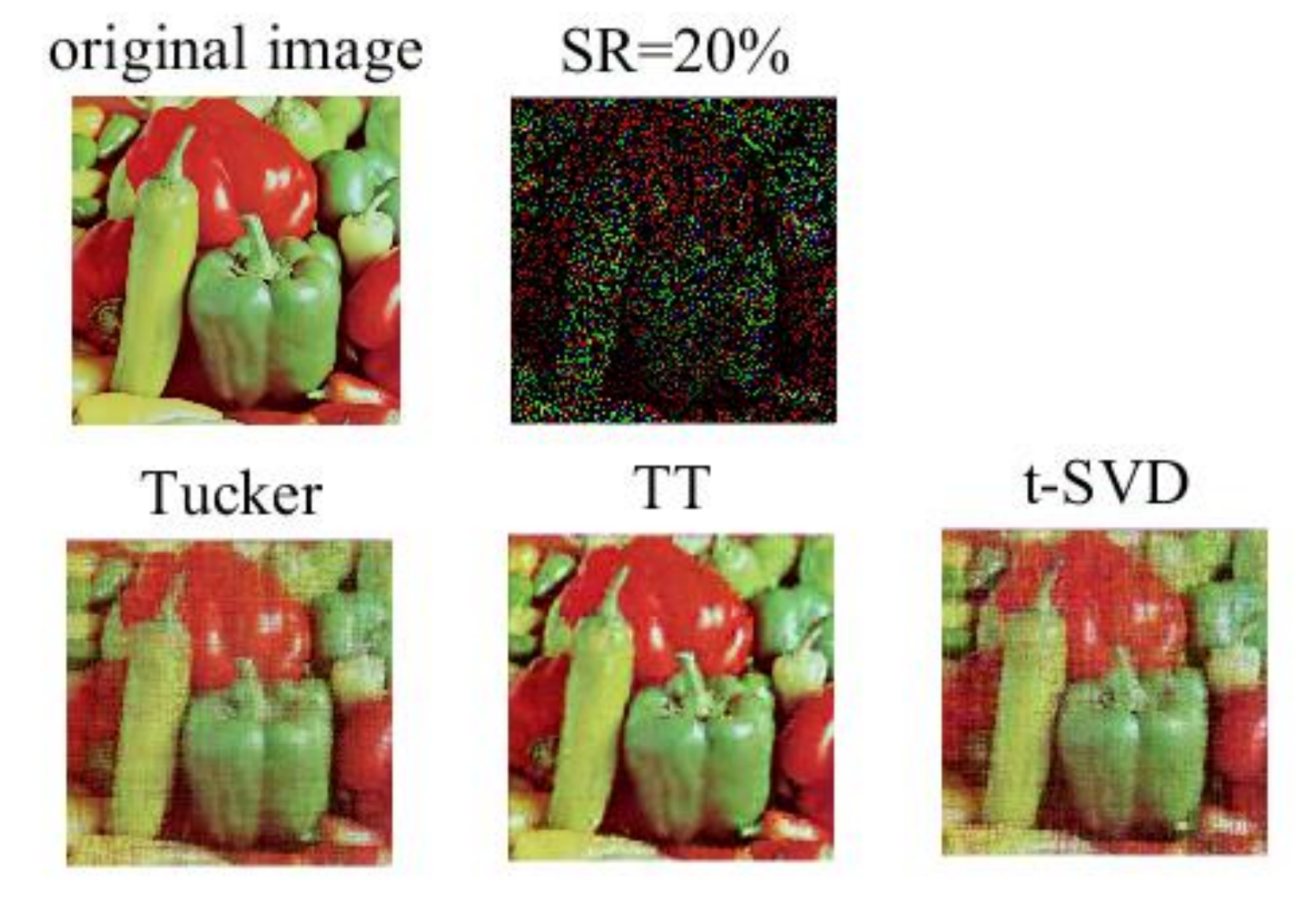}
	\caption{The recovered results with different rank minimization.}
	\label{MinRank}
\end{figure}

\begin{figure}[t]
	\centering
	\includegraphics[width=300pt, keepaspectratio]{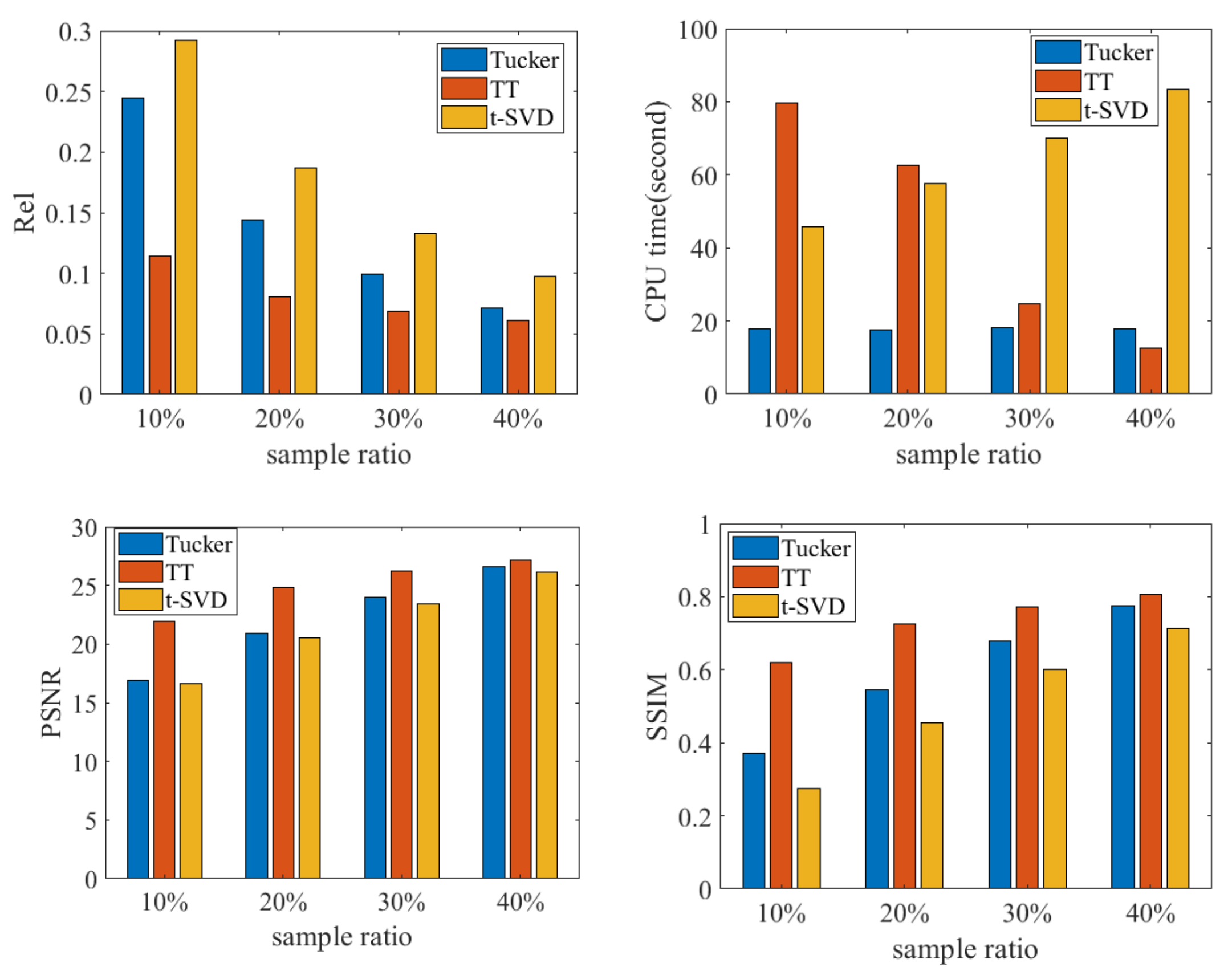}
	\caption{The performance of different decomposition method with rank minimization in terms of Rel, CPU time, PSNR and SSIM.}
	\label{Relrank}
\end{figure}

Fig. \ref{MinRank} are recovered results with different rank minimization under the $\text{SR}=20\% $. And Fig. \ref{Relrank} shows the performance of different rank minimization with respect to Rel, PSNR and SSIM under the SR changes from 10\% to 40\%.

\section{Conclusions}
\label{sec:6}
In this review, we provide clear guidelines of tensor completion for researchers in processing color images and videos. First, we divided tensor completion algorithms into two groups according to predefined rank and rank minimization, and then for each group, we introduce a variety of optimization models in details with respect to different tensor decomposition. Additionally, we conduct a group of experiments using the factorization based approaches and rank minimization models. And we find that the rank minimization models perform well in terms of the accuracy and cost time.

\section{Discussion and Future Work}
\label{sec:7}
We introduce a variety of tensor completion methods based on different tensor decomposition. Traditional tensor decompositions, such as CP decomposition and TucIker decomposition, have been studied many times. The CP decomposition can reduce the storage complexity from $I^{D}$ to $IDS$ where $D$ is the dimension of a tensor, $I$ is the size corresponding to dimension and $S$ is the CP rank, but the minimization of tensor nuclear norm for CP rank is an NP-hard problem. Hence, the existing low CP rank tensor completion methods often iteratively update the factor matrices with a predefined rank. The Tucker decomposition can reduce the storage complexity from $I^{D}$ to $IDS+S^{D}$, where  $D$ is the dimension of a tensor, $I$ is the size corresponding to dimension， and we assume $S_{1}=S_{2}=\dots=S_{N}=S$, $\operatorname{rank}_{\text{Tucker}}=(S_{1},S_{2},\dots,S_{N})$, which indicts that the storage complexity grows exponentially with the dimensional increasing. In addition, t-SVD is a good extending of matrix, but it only suits for 3rd-order tensor. The simple tensor network decomposition including tensor train decomposition and tensor ring decomposition has shown its performance with traditional tensor composition, but the related work is little.

With the increasing of the data size, the computation cost of traditional tensor decomposition such as Tucker decomposition and CP decomposition is more expensive. In contrast, tensor network decomposition can replace it. Besides, more efficient and fast algorithms need to be proposed to deal with the big data. Furthermore, applying the tensor space into a higher-order space is one of the important work in the future.
\section{Acknowledgment}
\label{sec: Acknowledgment}
The authors would like to  thank Jiani Liu, Jiayi Fang and Huyan Huang for the helpful discussions. This research is supported by National Natural Science Foundation of China (NSFC, No. 61602091, No. 61571102) and the Fundamental Research Funds for the Central Universities (No. ZYGX2016J199, No. ZYGX2014Z003).

\end{document}